\documentclass[reqno]{amsart}

\usepackage{pstricks}
\usepackage{verbatim}
\usepackage{amssymb}
\usepackage{mathrsfs}
\usepackage{dsfont}
\usepackage{amsmath}
\usepackage{amsthm}

\AtBeginDocument{%
}


\title{Spherical Harmonic Analysis on Affine Buildings}

\numberwithin{equation}{section}

\newtheorem{lem}{Lemma}[section]
\newtheorem{thm}[lem]{Theorem}
\newtheorem{cor}[lem]{Corollary}
\newtheorem{prop}[lem]{Proposition}

\theoremstyle{definition}
\newtheorem{defn}[lem]{Definition}

\newtheorem{rem}[lem]{Remark}

\newtheorem*{ack}{Acknowledgements}

\newcommand{\taua}{\tau_{\alpha\vphantom{/2}}}
\newcommand{\tauah}{\tau_{\alpha/2}^{-1/2}}
\newcommand{\Hom}{\mathrm{Hom}}

\renewcommand{\S}{\Sigma}
\renewcommand{\t}{\theta}
\newcommand{\la}{\lambda}

\renewcommand{\o}{\omega}
\renewcommand{\O}{\Omega}
\newcommand{\s}{\sigma}
\newcommand{\ca}{\mathcal{A}}
\newcommand{\cb}{\mathcal{B}}
\newcommand{\cc}{\mathcal{C}}

\newcommand{\ch}{\mathcal{H}}

\newcommand{\cs}{\mathcal{S}}
\newcommand{\ct}{\mathcal{T}}

\newcommand{\cn}{\mathcal{N}}

\newcommand{\bc}{\mathbb{C}}
\newcommand{\bn}{\mathbb{N}}
\newcommand{\br}{\mathbb{R}}
\newcommand{\bt}{\mathbb{T}}
\newcommand{\bz}{\mathbb{Z}}

\newcommand{\bq}{\mathbb{Q}}

\newcommand{\sca}{\mathscr{A}}

\newcommand{\lan}{\langle}
\newcommand{\ran}{\rangle}

\newcommand{\Aut}{\mathrm{Aut}}
\newcommand{\Auttr}{\mathrm{Aut_{tr}}}

\newcommand{\scx}{\mathscr{X}}

\newcommand{\ts}{\textsection}

\newcommand{\cha}{\alpha^{\vee}}
\newcommand{\chR}{R^{\vee}}

\newcommand{\scc}{\mathscr{C}}

\DeclareMathAlphabet{\mathpzc}{OT1}{pzc}{m}{it}

\newcommand{\e}{\epsilon}
\newcommand{\bu}{\mathbb{U}}

\begin{document}

\author{James Parkinson}



\dedicatory{\upshape
  School of Mathematics and Statistics\\
  University of Sydney, NSW 2006\\
  Australia\\[3pt]
  \texttt{jamesp@maths.usyd.edu.au}\\[8pt]
  \today}

\begin{abstract} Let $\scx$ be a locally finite regular affine building
with root system~$R$. There is a commutative algebra~$\sca$
spanned by averaging operators~$A_{\la}$, $\la\in P^+$, acting on
the space of all functions $f:V_P\to\bc$, where $V_P$ is in most
cases the set of all special vertices of $\scx$, and $P^+$ is a
set of dominant coweights of~$R$. This algebra is studied in \cite{C3} 
and~\cite{C2} for $\tilde{A}_n$ buildings,
and the general case is treated in~\cite{p}.

In this paper we show that all algebra homomorphisms
$h:\sca\to\bc$ may be expressed in terms of the \textit{Macdonald
spherical functions}. We also provide a second formula for these
homomorphisms in terms of an integral over the \textit{boundary}
of $\scx$. We may regard $\sca$ as a subalgebra of the
$C^*$-algebra of bounded linear operators on $\ell^2(V_P)$, and we
write $\sca_2$ for the closure of $\sca$ in this algebra. We study
the Gelfand map $\sca_2\to\scc(M_2)$, where
$M_2=\Hom(\sca_2,\bc)$, and we compute $M_2$ and the
\textit{Plancherel measure} of $\sca_2$. We also compute the
$\ell^2$-operator norms of the operators $A_{\la}$, $\la\in P^+$,
in terms of the Macdonald spherical functions.
\end{abstract}

\maketitle

\section*{Introduction}

In \cite{p} we showed that under the weak hypothesis of
regularity, to each affine building $\scx$ of irreducible type
there is naturally associated a commutative algebra $\sca$ spanned
by \textit{vertex set averaging operators} $A_{\la}$, $\la\in
P^+$, acting on the space of all functions $f:V_P\to\bc$. Here
$P^+$ is the set of dominant coweights of an appropriate root
system $R$, and $V_P$ is in most cases the set of all
\textit{special vertices} of $\scx$.

Let $\bc[P]$ denote the group algebra of the coweight lattice $P$
of $R$, with the group operation written multiplicatively. Thus
elements of $\bc[P]$ are linear combinations of the formal
exponentials $\{x^{\la}\}_{\la\in P}$, and
$x^{\la}x^{\mu}=x^{\la+\mu}=x^{\mu}x^{\la}$ for all $\la,\mu\in
P$. Let $W_0$ denote the \textit{Weyl group} of $R$. Thus $W_0$
acts on $P$, and hence on $\bc[P]$ by linearly extending
$wx^{\la}=x^{w\la}$. Let $\bc[P]^{W_0}=\{f\in\bc[P]\mid
wf=f\textrm{ for all $w\in W_0$}\}$. In \cite[Theorem~6.16]{p} we
showed that $\sca\cong\bc[P]^{W_0}$, with the isomorphism given by
$A_{\la}\mapsto P_{\la}(x)$, where the elements
$P_{\la}(x)\in\bc[P]^{W_0}$, $\la\in P^+$, are the
\textit{Macdonald spherical functions}.

In this paper we study the algebra homomorphisms $h:\sca\to\bc$.
The isomorphism $\sca\cong\bc[P]^{W_0}$ gives one formula for
these homomorphisms in terms of the Macdonald spherical functions.
The main part of this paper is devoted to a rather different line
of reasoning which furnishes a second formula for the algebra
homomorphisms $h:\sca\to\bc$ in terms of an integral over the
\textit{boundary} of $\scx$ (that is, the \textit{spherical
building at infinity}). This formula is an analogue of the formula
in \cite[Proposition~3.3.1]{macsph}, which expresses the
\textit{zonal spherical functions on a group $G$ of $p$-adic type}
as an integral over $K$, where $K$ is a certain compact subgroup
of $G$. For example, $G=SL(n+1,F)$, where $F$ is a $p$-adic field,
and $K=SL(n+1,\mathcal{O})$, where $\mathcal{O}$ is the ring of
integers in $F$. See below for more details.

The algebra $\sca$ may be regarded as a subalgebra of the
$C^*$-algebra $\mathscr{L}(\ell^2(V_P))$ of bounded linear
operators on $\ell^2(V_P)$, and the closure $\sca_2$ of $\sca$ in
$\mathscr{L}(\ell^2(V_P))$ is a commutative $C^*$-algebra. We
study the Gelfand map $\sca_2\to\scc(M_2)$, where
$M_2=\Hom(\sca_2,\bc)$. We compute $M_2$ and the Plancherel
measure $\pi$ on $M_2$ corresponding to the natural state
$A\mapsto (A\delta_o)(o)$ on $\sca_2$ (here $o\in V_P$ is any
fixed vertex, and $\delta_o(x)=1$ if $x=o$ and $0$ otherwise).
These results will be used both in this paper (for example, in the
important Theorem~\ref{mainresult} where we show the equality of
our two expressions for the homomorphisms $h:\sca\to\bc$) and in
later work, where we will study \textit{isotropic random walks} on
affine buildings. In Theorem~\ref{spectralradii} we compute the
$\ell^2$-operator norms $\|A_{\la}\|$, $\la\in P^+$, in terms of
the Macdonald spherical functions.

It is important to notice that throughout this work we do not
assume the existence of any groups acting in particular ways on
the building. We only use the building axioms, which puts our
results into a general setting.

Let us briefly describe the connection between our work and that
in~\cite{macsph}. Let~$G$ be a group of $p$-adic type with a
maximal compact subgroup~$K$, as in~\cite[\ts2.7]{macsph}.
Associated to~$G$ there is an irreducible (but not necessarily
reduced) root system~$R$, as in~\cite[Chapter~II]{macsph}. As
mentioned above, the typical example here is \mbox{$G=SL(n+1,F)$}
and \mbox{$K=SL(n+1,\mathcal{O})$}, where $F$ is a $p$-adic field
and $\mathcal{O}$ is the ring of integers in~$F$. In this
case~$R$ is a root system of type~$A_n$.

A function $f:G\to\bc$ is called \textit{bi-$K$-invariant} if
$f(gk)=f(kg)=f(g)$ for all $g\in G$ and $k\in K$. Let
$\mathscr{L}(G,K)$ be the space of continuous, compactly supported
bi-$K$-invariant functions on~$G$. In
\cite[Theorem~3.3.6]{macsph}, Macdonald shows that
$\mathscr{L}(G,K)\cong \bc[Q]^{W_0}$ (the subalgebra of
$W_0$-invariant elements of the group algebra of the coroot
lattice~$Q$ of~$R$), and thus $\mathscr{L}(G,K)$ is a commutative
convolution algebra.

A function $\phi:G\to\bc$ is called a \textit{zonal spherical
function relative to~$K$} if
\begin{enumerate}
\item[(i)] $\phi(1)=1$,
\item[(ii)] $\phi$ is bi-$K$-invariant and continuous, and
\item[(iii)] $f*\phi=\la_f\phi$ for all $f\in\mathscr{L}(G,K)$,
where $\la_f$ is a scalar
\end{enumerate}
(see \cite[Proposition~1.2.5]{macsph}). In
\cite[Proposition~3.3.1]{macsph} Macdonald gives a formula for
the zonal spherical functions in terms of an integral over~$K$,
and in \cite[Theorem~4.1.2]{macsph} he uses this integral formula
to obtain a second `summation' formula for the spherical functions
in terms of a sum over~$W_0$ (the Weyl group of~$R$) of rational
functions.

The group $G$ acts strongly transitively on its Bruhat-Tits
building~$\scx$ (see \cite[Lemma~2.4.4]{macsph}), which is locally
finite, regular and affine, although Macdonald makes little use of
this building. The contents of this paper (and \cite{p}) lead us
to the conclusion that, rather than playing a relatively minor
role, the building theoretic elements alone determine the nature
of the algebra $\mathscr{L}(G,K)$ and the zonal spherical
functions.

Specifically, the algebra $\mathscr{L}(G,K)$ is isomorphic to a
subalgebra $\sca_0$ of $\sca$, spanned by the operators
\mbox{$\{A_{\la}\mid\la\in Q\cap P^+\}$}. The reason that this
smaller algebra occurs here is that Macdonald supposes that his
groups of $p$-adic type are simply connected, and so for him the
coweight lattice has to be replaced by the coroot lattice, and
thus the set $V_P$ is replaced by the smaller set $V_Q$ consisting
of all vertices of one special type. See
\cite[Propositions~2.4~and~2.5]{cm} for details in the~$R=A_2$
case. In~\cite[Theorem~6.16]{p} we showed that
$\sca\cong\bc[P]^{W_0}$ (and indeed one can easily see that
$\sca_0\cong\bc[Q]^{W_0}$), thus proving the analogue
of~\cite[Theorem~3.3.6]{macsph}.

The zonal spherical functions on $G$ correspond to the spherical
functions on~$V_Q$ (see Definition~\ref{zsf}), which in turn
correspond to the algebra homomorphisms \mbox{$h:\sca_0\to\bc$}
(see Proposition~\ref{conn}). Our analogue of Macdonald's integral
formula is Theorem~\ref{main} (see also Corollary~\ref{2ndform}),
and our analogue of his summation formula is~(\ref{macsph}). Both
of our formulae are proved in the context of the larger algebra
$\sca$ (in fact we will have little to say regarding the smaller
algebra~$\sca_0$ in what follows). Since there are locally finite
regular affine buildings that are not the Bruhat-Tits buildings
of any group, our building approach puts the results of
\cite{macsph} into a more general setting.

\begin{ack} The author would like to thank Donald Cartwright for many
helpful suggestions throughout the preparation of this paper.
\end{ack}

\section{Notations and Definitions}\label{highestroot}

\subsection{Root Systems} Root systems play a significant role in this
work. We fix the following notations and conventions, generally
following \cite[Chapter~VI]{bourbaki}.

Let $R$ be an irreducible, but not necessarily reduced, root
system in a real vector space $E$ with inner product
$\lan\cdot,\cdot\ran$. The \textit{rank} of $R$ is $n$, the
dimension of $E$. Let $I_0=\{1,2,\ldots,n\}$ and
$I=\{0,1,2,\ldots,n\}$. Let $B=\{\alpha_i\mid i\in I_0\}$ be a
fixed base of $R$, and write $R^+$ for the set of positive roots
(relative to $B$). For $\alpha\in E\backslash\{0\}$, write
$\alpha^{\vee}=\frac{2\alpha}{\lan\alpha,\alpha\ran}$, and let
$R^{\vee}=\{\alpha^{\vee}\mid\alpha\in R\}$ be the \textit{dual}
root system of~$R$. Since $R$ is irreducible there is a unique
\textit{highest root} $\tilde{\alpha}$, and we define numbers
$\{m_{i}\}_{i\in I}$ by $m_0=1$, and $\tilde{\alpha}=\sum_{i\in
I_0}m_i\alpha_i$. For each $i\in I_0$ define $\la_i\in E$ by $\lan
\la_i,\alpha_j\ran=\delta_{i,j}$. The elements $\{\la_i\}_{i\in
I_0}$ are called the \textit{fundamental coweights of $R$}. The
\textit{coweight lattice of $R$} is $P=\sum_{i\in I_0}\bz\la_i$,
and elements $\la\in P$ are called \textit{coweights of $R$}. A
coweight $\la\in P$ is said to be \textit{dominant} if
$\lan\la,\alpha_i\ran\geq0$ for all $i\in I_0$, and we write $P^+$
for the set of all dominant coweights. Finally, let
$Q=\sum_{\alpha\in R}\bz\alpha^{\vee}$ be the \textit{coroot
lattice of $R$}, and let $Q^+=\sum_{\alpha\in
R^+}\bn\alpha^{\vee}$, where $\bn=\{0,1,2\ldots\}$. Note that
$Q\subseteq P$.

For those readers not so familiar with the non-reduced root
systems, we make the following comments. For each $n\geq1$ there
is exactly one irreducible non-reduced root system (up to
isomorphism) of rank $n$, denoted by $BC_n$. To describe the root
system $BC_n$, we may take $E=\br^n$ with the usual inner
product, and let $R$ consist of the vectors $\pm e_i,\pm 2e_i$ and
$\pm e_j\pm e_k$ for $1\leq i\leq n$ and $1\leq j<k\leq n$. Let
$\alpha_j=e_j-e_{j+1}$ for $1\leq j<n$ and $\alpha_n=e_n$. Then
$B=\{\alpha_j\}_{j=1}^{n}$ is a base of $R$, and $R^+$ consists
of the vectors $e_i,2e_i$ and $e_j\pm e_k$ for $1\leq i\leq n$
and $1\leq j<k\leq n$. The fundamental coweights are
$\la_i=e_1+\cdots+e_i$ for each $1\leq i\leq n$. Note that
$\chR=R$ and $Q=P$. The subsystem $R_1=\{\alpha\in R\mid
2\alpha\notin R\}$ is a root system of type $C_n$, and the
subsystem $R_2=\{\alpha\in R\mid \frac{1}{2}\alpha\notin R\}$ is
a root system of type $B_n$ (with the convention that
$C_1=B_1=A_1$).

\subsection{Weyl Groups} For each $\alpha\in R$ and $k\in\bz$ let
$H_{\alpha;k}=\{x\in E\mid\lan x,\alpha\ran=k\}$. We call these
sets \textit{affine hyperplanes}, or simply \textit{hyperplanes}.
For each $\alpha\in R$ and $k\in \bz$ let $s_{\alpha;k}$ denote
the orthogonal reflection in $H_{\alpha;k}$. Thus
$s_{\alpha;k}(x)=x-(\lan x,\alpha\ran-k)\alpha^{\vee}$ for all
$x\in E$. Write $s_{\alpha}$ in place of $s_{\alpha;0}$, $s_i$ in
place of $s_{\alpha_i}$ (for $i\in I_0$), and let
$s_0=s_{\tilde{\alpha};1}$. Let $W_0=W_0(R)$ be the \textit{Weyl
group of $R$}, and let $W=W(R)$ be the \textit{affine Weyl group
of $R$}. Thus $W_0$ is the subgroup of $\mathrm{GL}(E)$ generated
by $S_0=\{s_i\}_{i\in I_0}$, and $W$ is the subgroup of
$\mathrm{Aff}(E)$ generated by $S=\{s_i\}_{i\in I}$. Both
$(W_0,S_0)$ and $(W,S)$ are Coxeter systems. For subsets $J\subset
I$ we write $W_{J}$ for the subgroup of $W$ generated by
$\{s_i\}_{i\in J}$. Thus $W_0=W_{I\backslash\{0\}}$. Given $w\in
W$, we define the \textit{length} $\ell(w)$ of $w$ to be smallest
$n\in\bn$ such that $w=s_{i_1}\ldots s_{i_n}$, with
$i_1,\ldots,i_n\in I$. We write $M=(m_{i,j})_{i,j=0}^n$ for the
\textit{Coxeter matrix} of $W$. Thus $m_{i,j}$ is the order of
$s_is_j$ in $W$.

\subsection{The Coxeter Complex}
The \textit{Coxeter complex} of $W$ is the \textit{labelled
simplicial complex} $\S_s$ whose vertex set is the (automatically
disjoint) union over $i\in I$ of the sets
$X_i=\{wW_{I\backslash\{i\}}\mid w\in W\}$. If $x\in X_i$ we say
that $x$ has \textit{type $i$}, and write $\tau(x)=i$. Simplices
of $\S_s$ are subsets of the \textit{maximal simplices} (or
\textit{chambers}), which are defined to be sets of the form
$\{wW_{I\backslash\{i\}}\mid i\in I\}$, $w\in W$. The subscript
`$s$' on $\S_s$ stands for \textit{`simplicial'}. It is also
possible to think of $\S_s$ as a \textit{chamber system} $\S_c$.
Let $\S_c$ be the set of all chambers of $\S_s$, and for
$C,C'\in\S_c$ and $i\in I$, declare $C\sim_i C'$ if either $C=C'$
or if all the vertices of $C$ and $C'$ are the same, except for
those of type~$i$.

There is a natural \textit{geometric realisation} $\S$ of $\S_s$
in the case when $R$ is irreducible (when $R$ is reducible the
following construction produces a \textit{polysimplicial
complex}). Let $\ch$ denote the family of hyperplanes
$H_{\alpha;k}$, $\alpha\in R$, $k\in\bz$, and define
\textit{chambers} of $\S$ to be open connected components of
$E\backslash\bigcup_{H\in \ch}H$. Since $R$ is irreducible, each
chamber is an open (geometric) simplex \cite[V, \ts3, No.9,
Proposition~8]{bourbaki}, and we call the extreme points of the
closure of chambers \textit{vertices} of $\S$. We write $V(\S)$
for the set of all vertices of $\S$.

The choice of the base $B$ gives a natural choice of a
\textit{fundamental chamber}
\begin{align}
\label{fundamentalalcove}C_0=\{x\in E\mid\lan
x,\alpha_i\ran>0\textrm{ for all $i\in I_0$ and }\lan
x,\tilde{\alpha}\ran<1\}\,.
\end{align}
The vertices of $C_0$ are the points $\{0\}\cup\{\la_i/m_i\}_{i\in
I_0}$ \cite[VI, \ts2, No.2]{bourbaki}. There is a natural
simplicial complex structure on $\S$ by taking maximal simplices
to be the vertex sets of chambers of $\S$, and taking simplices
to be subsets of maximal simplices. Thus, following the argument
in \cite[Lemma~1.5]{donaldintro}, we define $\tau:V(\S)\to I$ to
be the unique labelling of $\S$ (as a simplicial complex) such
that $\tau(0)=0$ and $\tau(\la_i/m_i)=i$. As a labelled simplicial
complex, $\S$ is isomorphic to $\S_s$.

We will henceforth simply write $\S$ for the Coxeter complex of
$W$. It will be clear from the context if one is to regard $\S$ as
a simplicial complex or as a geometric realisation in $E$.

\subsection{The `Good' Vertices and Sectors of $\S$} At this stage there is really no
difference between the $C_n$ and $BC_n$ cases. That is,
$\S_s(C_n)\equiv \S_s(BC_n)$. We now introduce some additional
structure that will differentiate between these two cases.

The set $P$ of coweights of $R$ is a subset of $V(\S)$, and we
call elements of $P$ the \textit{good vertices of $\S$}. Note that
$P(C_n)\neq P(BC_n)$. When $R$ is reduced, $P$ is the set of more
familiar \textit{special vertices} of $\S$ \cite[VI, \ts2, No.2,
Proposition~3]{bourbaki}.

We write $I_{P}=\{\tau(\la)\mid\la\in P\}\subseteq I$. We have
$I_P=\{i\in I\mid m_i=1\}$, which shows that $0\in I_P$ for all
root systems, and that $I_P=\{0\}$ if $R$ is non-reduced
\cite[Lemma~4.3]{p}. This also shows that in the non-reduced case
there are `special' vertices which are not `good' vertices.

We define the \textit{fundamental sector} of $\S$ to be the open
simplicial cone
\begin{align}
\label{fundamentalsector}\cs_0=\{x\in E\mid\lan
x,\alpha_i\ran>0\textrm{ for all $i\in I_0$}\}.
\end{align}
The \textit{sectors} of $\S$ are then the sets $\la+w(\cs_0)$,
where $w\in W_0$ and $\la\in P$.

\subsection{The Extended Affine Weyl Group} Let $\tilde{W}=\tilde{W}(R)$ be the \textit{extended affine Weyl
group of $R$}. That is, $\tilde{W}=W_0\ltimes P$. For $w\in
\tilde{W}$ we define
$$\ch(\tilde{w})=\{H\in\ch\mid H\textrm{ separates $C_0$ and
$\tilde{w}C_0$}\},$$ and we define the \textit{length of
$\tilde{w}$} to be $\ell(\tilde{w})=|\ch(\tilde{w})|$. For
$\la\in P$, we write $t_{\la}\in \tilde{W}$ for the translation
$t_{\la}(x)=\la+x$ for all $x\in E$.

Let $G=\{g\in\tilde{W}\mid\ell(g)=0\}$; this is the stabiliser of $C_0$ in
$\tilde{W}$. We have
$\tilde{W}\cong W\rtimes G$ \cite[VI, \ts2, No.3]{bourbaki}, and
furthermore $G\cong P/Q$, and so $G$ is a finite abelian group.
For each $\la\in P$, let $W_{0\la}$ denote the stabiliser of
$\la$ in $W_0$. Let $w_0$ and $w_{0\la_i}$ denote the longest
elements of $W_0$ and $W_{0\la_i}$ respectively. Then
\begin{align}
\label{G}G=\{g_i\mid i\in I_P\}
\end{align}
where $g_0=1$ and $g_i=t_{\la_i}w_{0\la_i}w_0$ for $i\in
I_P\backslash\{0\}$ \cite[(4.5)]{p}.

\subsection{Automorphisms} An \textit{automorphism of $\S$} is a
bijection $\psi$ of $E$ which maps chambers, and only chambers, to
chambers, with the property that chambers $C$ and $C'$ are
adjacent if and only if $\psi(C)$ is adjacent to~$\psi(C')$. We
write $\Aut(\S)$ for the automorphism group of $\S$.

Let $D$ be the Coxeter diagram of $W$. An \textit{automorphism of
$D$} is a permutation $\s$ of $I$ such that
$m_{\s(i),\s(j)}=m_{i,j}$ for all $i,j\in I$. We write $\Aut(D)$
for the automorphism group of $D$.

For each $\psi\in\Aut(\S)$ there exists a $\s\in\Aut(D)$ such
that $(\tau\circ\psi)(v)=(\s\circ\tau)(v)$ for all $v\in V(\S)$,
and if $C\sim_i C'$ then $\psi(C)\sim_{\s(i)}\psi(C')$
\cite[Proposition~4.6]{p}. An automorphism $\psi\in\Aut(\S)$ is
called \textit{type preserving} if the associated $\s\in\Aut(D)$
is the identity. By \cite[Lemma~2.2]{ronan} $\psi\in\Aut(\S)$ is
type preserving if and only if~$\psi\in W$.

The group $\tilde{W}$ acts on $\S$ and so embeds as a subgroup of
$\Aut(\S)$. The automorphisms $\tilde{w}\in\tilde{W}$, and the
associated automorphisms $\s\in\Aut(D)$, are called \textit{type
rotating}. We write $\Auttr(D)$ from the subgroup of $\Aut(D)$
consisting of all type rotating automorphisms of $D$. The group
$\Auttr(D)$ may be described as follows. For each $i\in I_P$, let
$\s_i$ be the automorphism of $D$ defined by the equation
$\s_i\circ\tau=\tau\circ g_i$. Since each $\tilde{w}\in\tilde{W}$
may be written as $\tilde{w}=wg_i$ where $w\in W$ and $i\in I_P$
\cite[VI, \ts2, No.3]{bourbaki}, we have $\Auttr(D)=\{\s_i\mid
i\in I_P\}$. The group $\Auttr(D)$ is abelian (since $G$ is), it
acts simply transitively on $I_P$, and $\s_i(0)=i$ for each $i\in
I_P$ (see \cite[Proposition~4.7(iii)]{p}).

\begin{rem} In the $\tilde{A}_n$ case, the qualification
\textit{type rotating} is very natural, for $\Auttr(D)$ consists
of the permutations $i\mapsto k+i\mod(n+1)$ for $k=0,\ldots,n$.
In the general case, the term \textit{type rotating} is less
natural.
\end{rem}

\subsection{Buildings and Regularity}\label{section:bar} Let $M$ be the Coxeter
matrix of $W$. Recall (\cite[p.76--77]{brown}) that a
\textit{building of type $M$} is a nonempty simplicial complex
$\scx$ which contains a family of subcomplexes called
\textit{apartments} such that
\begin{enumerate}
\item[$\mathrm{(i)}$] each apartment is isomorphic to the (simplicial) Coxeter
complex of $W$,
\item[$\mathrm{(ii)}$] given any two chambers of $\scx$ there is an
apartment containing both, and
\item[$\mathrm{(iii)}$] given any
two apartments $\ca$ and $\ca'$ that contain a common chamber,
there exists an isomorphism $\psi:\ca\to\ca'$ fixing $\ca\cap\ca'$
pointwise.
\end{enumerate}
We call $\scx$ \textit{affine} when $W$ is (as we assume
throughout).

It is an easy consequence of this definition that $\scx$ is a
labellable simplicial complex (with set of \textit{types}~$I$).
All of the isomorphisms in the above definition may be taken to be
label preserving (this ensures that the labellings of $\scx$ and
$\S$ are \textit{compatible}). The isomorphism in (iii) is then
unique~\cite[p.77]{brown}. Recall (\cite[p.30]{brown}) that the
\textit{type} of a simplex $\s$ is the subset of $I$ consisting of
the labels of the vertices of $\s$. If $\s$ has type $J\subset I$,
then the \textit{cotype} of $\s$ is $I\backslash J$.

We write $V$ for the vertex set of $\scx$, and $\cc$ for the set
of all chambers of $\scx$. As usual, we declare chambers
$c,d\in\cc$ to be \textit{$i$-adjacent}, and write $c\sim_i d$, if
and only if either $c=d$, or if all the vertices of $c$ and $d$
are the same, except for those of type~$i$. By a
\textit{gallery} of type $i_1\cdots i_n\in I^*$
(here $I^*$ denotes the free monoid on~$I$) in $\scx$ we mean a
sequence $c_0,\ldots,c_n$ of chambers such that
$c_{k-1}\sim_{i_k}c_k$ and $c_{k-1}\neq c_k$ for $1\leq k\leq n$.
If we remove the condition that $c_{k-1}\neq c_k$ for each $1\leq
k\leq n$, we call the sequence $c_0,\ldots,c_n$ a
\textit{pre-gallery}.

We say that $\scx$ is \textit{locally finite} if both $|I|<\infty$
and $|\{d\in\cc\mid d\sim_i c\}|<\infty$ for all $i\in I$ and
$c\in\cc$. We call $\scx$ \textit{regular} if for each $i\in I$
the cardinality $|\{d\in\cc\mid d\sim_i c\}|$ is independent of
$c\in\cc$. In this case we define numbers $\{q_i\}_{i\in I}$ by
$q_i+1=|\{d\in\cc\mid d\sim_i c\}|$. The numbers $q_i$, $i\in I$,
are called the \textit{parameters} of the building. These
parameters satisfy $q_j=q_i$ if $s_j=ws_iw^{-1}$ for some $w\in W$
(see \cite[Corollary~2.2]{p}). If $w=s_{i_1}\cdots s_{i_m}\in W$
is a reduced expression we define $q_w=q_{i_1}\cdots q_{i_m}$,
which is independent of the particular reduced expression for $w$
(see \cite[Proposition~2.1(i)]{p} or \cite[IV, \ts1, No.5,
Proposition~5]{bourbaki}).

To each locally finite regular affine building (of irreducible
type) we associate an irreducible root system~$R$ (depending on
the parameter system of the building) as follows:
\begin{enumerate}
\item[(i)] If $\scx$ is a regular building of type $\tilde{X}_n$, where
$X=A$ and $n\geq2$, or $X=B$ and $n\geq3$, or $X=D$ and $n\geq4$,
or $X=E$ and $n=6,7$ or $8$, or $X=F$ and $n=4$, or $X=G$ and
$n=2$, then we take $R=X_n$.
\item[(ii)] If $\scx$ is a regular $\tilde{A}_1$ building with $q_0\neq q_1$, then we
take $R=BC_1$.
\item[(iii)] If $\scx$ is a regular $\tilde{C}_n$ building with $n\geq2$ and $q_0\neq q_n$, then we take
$R=BC_n$.
\item[(iv)] If $\scx$ is a regular $\tilde{A}_1$ building with
$q_0=q_1$, then we take $R=A_1$.
\item[(v)] If $\scx$ is a regular $\tilde{C}_n$ building with
$n\geq2$ and $q_0=q_n$, then we take $R=C_n$.
\end{enumerate}

The choices above are made to ensure that $\Auttr(D)$ preserves
the parameter system of $\scx$; that is, for $\s\in\Auttr(D)$,
$q_{\s(i)}=q_i$ for all $i\in I$. Thus, for example,~(iii) above
is motivated by the general parameter system of a $\tilde{C}_n$
building:\newline
\begin{figure}[ht]
 \begin{center}
 \psset{xunit= 0.9 cm,yunit= 0.9 cm}
 \psset{origin={0,0}}
\vspace{0.3cm}

\rput(-3,0.3){$q_0$}\pscircle*(-3,0){2pt}
\rput(-2,0.3){$q_1$}\pscircle*(-2,0){2pt} \psline(-3,0)(-2,0)
\rput(-2.5,0.3){$4$} \rput(-1,0.3){$q_1$}\pscircle*(-1,0){2pt}
\psline(-2,0)(-1,0) \pscircle*(0.4,0){1.5pt}
\pscircle*(-0.4,0){1.5pt} \pscircle*(0,0){1.5pt}
\rput(1,0.3){$q_1$}\pscircle*(1,0){2pt}
\rput(2,0.3){$q_1$}\pscircle*(2,0){2pt} \psline(1,0)(3,0)
\rput(3,0.3){$q_n$}\pscircle*(3,0){2pt} \rput(2.5,0.3){$4$}
\end{center}
\end{figure}

\noindent(see \cite[Appendix]{p}). If we take $R=C_n$, then
$\Auttr(D)=\{1,\s\}$, where $\s(0)=n$, and the condition
$q_{\s(0)}=q_0$ fails. However, if $R=BC_n$ then
$\Auttr(D)=\{1\}$.

If $R$ is the root system associated to $\scx$, we often say that
$\scx$ is an \textit{affine building of type~$R$}. Let $W_0$,
$W$, $\tilde{W}$, $P$, $Q$, and so on be as defined earlier.

We extend the definition of $q_w$ made above to $\tilde{W}$ by
setting $q_{\tilde{w}}=q_w$ if $\tilde{w}=wg$, where $w\in W$ and
$g\in G$. We note that if $\ell(uv)=\ell(u)+\ell(v)$ then
$q_{uv}=q_uq_v$. In particular, if $\la,\mu,\nu\in P^+$, then by
\cite[(2.4.1)]{m}
\begin{align}\label{ns2}
q_{t_{\la+\mu}}=q_{t_{\la}}q_{t_{\mu}},\quad\textrm{and so}\quad
q_{t_{\nu}}=\prod_{i=1}^n q_{t_{\la_i}}^{\lan\nu,\alpha_i\ran}.
\end{align}

\begin{defn}\label{goodbuilding} Let $\scx$ be an affine building of type~$R$ with
vertex set~$V$, and let~$\S=\S(R)$. Let $V_{\mathrm{sp}}(\S)$
denote the set of all special vertices of~$\S$, and let
$I_{\mathrm{sp}}=\{\tau(\la)\mid \la\in V_{\mathrm{sp}}(\S)\}$. Then:
\begin{enumerate}
\item[(i)] A vertex $x\in V$ is said to be \textit{special} if
$\tau(x)\in I_{\mathrm{sp}}$. We write $V_{\mathrm{sp}}$ for the
set of all special vertices of~$\scx$.
\item[(ii)] A vertex $x\in V$ is said to be \textit{good} if
$\tau(x)\in I_P$.
We write $V_P$ for the set of all good vertices of~$\scx$.
\end{enumerate}
\end{defn}

Clearly $V_P\subseteq V_{\mathrm{sp}}$. In fact if $R$ is reduced
then $V_P=V_{\mathrm{sp}}$. If $R$ is non-reduced (so $R$ is of
type $BC_n$ for some $n\geq1$), then $V_P$ is the set of all
type~$0$ vertices of~$\scx$, whereas $V_{\mathrm{sp}}$ is the set
of all type~$0$ and type~$n$ vertices of~$\scx$.

It is clear that a vertex $x\in V$ is good if and only if there exists an apartment~$\ca$ containing~$x$ and a type
preserving isomorphism~$\psi:\ca\to\S$ such that $\psi(x)\in P$.

\begin{defn} Let $\ca$ be an apartment of $\scx$. An isomorphism
$\psi:\ca\to\S$ is called \textit{type rotating} if and only if
it is of the form $\psi=w\circ\psi_0$, where $\psi_0:\ca\to\S$ is
a type preserving isomorphism, and $w\in \tilde{W}$. An
isomorphism $\varphi:\S\to\ca$ is called type rotating if and
only if $\varphi^{-1}:\ca\to\S$ is type rotating.
\end{defn}

\subsection{The Algebra $\sca$}\label{uu} Let $\scx$ be an affine building of type $R$ with parameter system
$\{q_i\}_{i\in I}$. Given $x\in V_P$ and $\la\in P^+$, let
$V_{\la}(x)$ be the set of all $y\in V_P$ such that there exists
an apartment $\ca$ containing $x$ and $y$, and a type rotating
isomorphism $\psi:\ca\to\S$ such that $\psi(x)=0$ and
$\psi(y)=\la$. In \cite[Proposition~5.6]{p} we showed that for
each $x\in V_P$, $\{V_{\la}(x)\}_{\la\in P^+}$ forms a partition
of~$V_P$.

\begin{rem} Let $V_Q$ denote the set of all type $0$ vertices of~$\scx$. If $\scx$ is the Bruhat-Tits building of a group~$G$
of $p$-adic type (see the introduction), then the sets $V_{\la}(x)$, $\la\in
Q\cap P^+$, are the orbits in $V_Q=G/G_x$
of the isotropy group~$G_x$.
\end{rem}

For $\la\in P^+$, let $\la^*=-w_0\la$, where $w_0$ is the longest
element of $W_0$. In \cite[Proposition~5.8]{p} we showed that
$\la^*\in P^+$ for all $\la\in P^+$, and that $y\in V_{\la}(x)$
if and only if $x\in V_{\la^*}(y)$. Note that $*$ is trivial
unless $w_0\neq-1$, that is, unless $R=A_n, D_{2n+1}$ or $E_6$ for
some $n\geq2$ (see \cite[Plates~I-IX]{bourbaki}).

Let $\la\in P^+$ and write $l=\tau(\la)$. Let $w_{\la}$ be the
unique minimal length double coset representative of
$W_{I\backslash\{0\}}t_{\la}'W_{I\backslash\{l\}}$, where
$t_{\la}'\in W$ is (uniquely) defined by $t_{\la}=t_{\la}'g$ for
some $g\in G$ (see \cite[\ts4.9]{p}). For finite subsets $U\subset
W$ we write
$$U(q)=\sum_{w\in U}q_w$$
for the \textit{Poincar\'{e} polynomial of $U$}.

In \cite[Theorem~5.15]{p} we showed that
$|V_{\la}(x)|=|V_{\la}(y)|$ for all $x,y\in V_P$ and $\la\in
P^+$, and we denote this common value by $N_{\la}$. We have
\begin{align}\label{firstformula}
N_{\la}=\frac{W_0(q)}{W_{0\la}(q)}q_{w_{\la}}=N_{\la^*},
\end{align} where $W_{0\la}=\{w\in W_0\mid w\la=\la\}$. We have
the following alternative formula.

\begin{prop}\label{betterf} Let $\la\in P^+$. Then
$$N_{\la}=\frac{W_0(q^{-1})}{W_{0\la}(q^{-1})}q_{t_{\la}},$$
where $W_0(q^{-1})=\sum_{w\in W_0}q_{w}^{-1}$, and similarly for
$W_{0\la}(q^{-1}).$
\end{prop}
\begin{proof} Let $l=\tau(\la)$, and let $w_0$ and $w_{0\la}$ be the longest elements of $W_0$ and
$W_{0\la}$, respectively. Since $t_{\la}=w_{\la}g_lw_0w_{0\la}$
and
$\ell(w_{\la}g_lw_0w_{0\la})=\ell(w_{\la}g_l)+\ell(w_0w_{0\la})$
(see \cite[Proposition~4.15]{p}) we have
$q_{t_{\la}}=q_{w_{\la}}q_{w_0w_{0\la}}$. By \cite[VI, \ts1, No.6,
Corollary~3 to Proposition~17]{bourbaki} we have
$\ell(w_0w)=\ell(w_0)-\ell(w)$ for all $w\in W_0$, and it follows
that $q_{w_0w_{0\la}}=q_{w_0}q_{w_{0\la}}^{-1}$. Thus
$$q_{w_{\la}}=q_{t_{\la}}q_{w_0w_{0\la}}^{-1}=q_{t_{\la}}q_{w_0}^{-1}q_{w_{0\la}}.$$

Furthermore, again using the fact that
$\ell(w_0w)=\ell(w_0)-\ell(w)$ for all $w\in W_0$, we have
$$W_0(q)=\sum_{w\in W_0}q_{w_0w}=q_{w_0}W_0(q^{-1}),$$
and similarly, since $W_{0\la}$ is a Coxeter group,
$W_{0\la}(q)=q_{w_{0\la}}W_{0\la}(q^{-1})$.

Putting all of this together, and using (\ref{firstformula}), the
result follows.
\end{proof}

We say that $\la\in P$ is \textit{strongly dominant} if
$\lan\la,\alpha_i\ran>0$ for all $i=1,\ldots,n$, and we write
$P^{++}$ for the set of all strongly dominant coweights. Note that
when $\la\in P^{++}$, $W_{0\la}=\{1\}$, and so by
Proposition~\ref{betterf}
\begin{align}
\label{ns}N_{\la}=W_0(q^{-1})q_{t_{\la}}\quad\textrm{for all
$\la\in P^{++}$.}
\end{align}

For each $\la\in P^+$ we define an operator $A_{\la}$, acting on
the space of all functions $f:V_P\to\bc$ by
\begin{align*}
(A_{\la}f)(x)=\frac{1}{N_{\la}}\sum_{y\in
V_{\la}(x)}f(y)\quad\textrm{for all $x\in V_P$.}
\end{align*}
Let $\sca$ be the linear span of $\{A_{\la}\}_{\la\in P^+}$ over
$\bc$. For each $\la,\mu\in P^+$ there exist numbers
$a_{\la,\mu;\nu}\in\bq^+$ such that
\begin{align*}
A_{\la}A_{\mu}=\sum_{\nu\in
P^+}a_{\la,\mu;\nu}A_{\nu}\quad\textrm{and }\sum_{\nu\in
P^+}a_{\la,\mu;\nu}=1
\end{align*}
(see \cite[Corollary~5.22]{p}), and so $\sca$ is an algebra. For
all $\la,\mu,\nu\in P^+$ we have
\begin{align}
\label{a}
a_{\la,\mu;\nu}=\frac{N_{\nu}}{N_{\la}N_{\mu}}|V_{\la}(x)\cap
V_{\mu^*}(y)|
\end{align}
where $x,y\in V_P$ are any vertices with $y\in V_{\nu}(x)$.

\subsection{The Isomorphism $A_{\la}\mapsto P_{\la}(x)$}

Let $\bc[P]$ be the group algebra (over $\bc$) of $P$, with the
group operation written multiplicatively. Thus elements of
$\bc[P]$ are linear combinations of the formal exponentials
$\{x^{\la}\}_{\la\in P}$. The group $W_0$ acts on $\bc[P]$ by
linearly extending $wx^{\la}=x^{w\la}$, and we write
$\bc[P]^{W_0}$ for the algebra consisting of all those $f\in
\bc[P]$ such that $wf=f$ for all $w\in W_0$.

Let $R_1=\{\alpha\in R\mid 2\alpha\notin R\}$, $R_2=\{\alpha\in
R\mid \frac{1}{2}{\alpha}\notin R\}$ and $R_3=R_1\cap R_2$ (so
$R_1=R_2=R_3=R$ if $R$ is reduced). For $\alpha\in R_2$, write
$q_{\alpha}=q_i$ if $|\alpha|=|\alpha_i|$, where $|\cdot|$ denotes
Euclidean root length (if $|\alpha|=|\alpha_i|$ then necessarily
$\alpha\in R_2$). Since $q_j=q_i$ whenever $s_j=ws_iw^{-1}$ for
some $w\in W$, it follows that $q_i=q_j$ whenever
$|\alpha_i|=|\alpha_j|$, and so the definition of $q_{\alpha}$ is
unambiguous.

Note that $R=R_3\cup(R_1\backslash R_3)\cup(R_2\backslash R_3)$
where the union is disjoint. Define a set of numbers
$\{\tau_{\alpha}\}_{\alpha\in R}$ by
\begin{align*}
\tau_{\alpha}=\begin{cases} q_{\alpha}&\textrm{if $\alpha\in
R_3$}\\ q_0&\textrm{if $\alpha\in R_1\backslash R_3$}\\
q_{\alpha}q_{0}^{-1}&\textrm{if $\alpha\in R_2\backslash
R_3$.}\end{cases}
\end{align*}
It is convenient to define $\tau_{\alpha}=1$ if $\alpha\notin R$.
Note that when $R$ is reduced, $\tau_{\alpha}=q_{\alpha}$, and
expressions like $\tau_{\alpha/2}$ and $\tau_{2\alpha}$ for
$\alpha\in R$ become~$1$ in subsequent formulae.

For $\la\in P^+$ define the \textit{Macdonald spherical function}
$P_{\la}(x)\in\bc[P]^{W_0}$ by
\begin{align}\label{eq:macsphfn}
P_{\la}(x)=\frac{q_{t_{\la}}^{-1/2}}{W_0(q^{-1})}\sum_{w\in
W_0}w\left(x^{\la}\prod_{\alpha\in R^+}\frac{1-\taua^{-1}\tauah
x^{-\cha}}{1-\tauah x^{-\cha}}\right).
\end{align}
By Proposition~\ref{equi}, this formula is equivalent to
\cite[(6.4)]{p}.

The main theorem of \cite{p} is the following.

\begin{thm}\cite[Theorem~6.16]{p}\label{me:fixref} The map $A_{\la}\mapsto P_{\la}(x)$ determines an algebra
isomomorphism, and so $\sca\cong\bc[P]^{W_0}$.
\end{thm}

\begin{rem} In \cite{p} we obtained the isomorphism
$\sca\to\bc[P]^{W_0}$ using the fact that $\bc[P]^{W_0}$ is the
center of an appropriately defined \textit{affine Hecke algebra}.
We will not need this fact here.
\end{rem}

\subsection{The Reducible Case}  When $\scx$ is a building of
reducible type, then $\scx$ has a natural description as a
\textit{polysimplicial complex}, and (assuming regularity) we can
associate a reducible root system $R$ to $\scx$ (bearing in mind
our conventions in Section~\ref{section:bar}). This leads to an
analogous definition of the algebra $\sca$ in this case.

It turns out that $\scx$ decomposes (essentially uniquely) into
the cartesian product of certain \textit{irreducible components}
$\{\scx_j\}_{j=1}^k$, each of which is an irreducible (regular)
building. Writing $\sca_j$ for the algebra of vertex set
averaging operators for $\scx_j$, it is not difficult to show
that $\sca\cong\sca_1\times\cdots\times\sca_k$, where $\times$ is
\textit{direct product} (see \cite[Remark~0.1]{p}). Thus
knowledge of the algebra homomorphisms $h_j:\sca_j\to\bc$ (using
the results of this paper) gives a complete description of the
algebra homomorphisms $h:\sca\to\bc$. Thus we can restrict our
attention to irreducible buildings throughout, without any loss of
generality.

\section{The Macdonald Formula for the Algebra Homomorphisms}

By Theorem~\ref{me:fixref}, the map $A_{\la}\mapsto P_{\la}(x)$
determines an isomorphism from $\sca$ onto $\bc[P]^{W_0}$. In this
section we will use this result to describe the \textit{Macdonald
formula} for the algebra homomorphisms $h:\sca\to\bc$.

For $u\in\Hom(P,\bc^{\times})$, write $u^{\la}$ in place of
$u(\la)$. Each $u\in\Hom(P,\bc^{\times})$ induces a homomorphism,
also called $u$, on $\bc[P]$, and all homomorphisms $\bc[P]\to\bc$
arise in this way. For $w\in W_0$ and $u\in\Hom(P,\bc^{\times})$
we write $wu\in\Hom(P,\bc^{\times})$ for the homomorphism
$(wu)^{\la}=u^{w\la}$. If $u\in\Hom(P,\bc^{\times})$, we write
$P_{\la}(u)$ in place of~$u(P_{\la}(x))$. Thus
\begin{align}\label{macsph}
\begin{aligned}
P_{\la}(u)&=\frac{q_{t_{\la}}^{-1/2}}{W_0(q^{-1})}\sum_{w\in
W_0}c(wu)u^{w\la},\quad\textrm{where}\\
c(u)&=\prod_{\alpha\in
R^+}\frac{1-\taua^{-1}\tauah u^{-\cha}}{1-\tauah u^{-\cha}},
\end{aligned}
\end{align}
provided, of course, that the denominators of the $c(wu)$
functions do not vanish. Since $P_{\la}(u)$ is a Laurent
polynomial, these \textit{singular cases} can be obtained from the
general formula by taking an appropriate limit (see
\cite[\ts4.6]{macsph}). Finally, for $u\in\Hom(P,\bc^{\times})$,
let $h_u:\sca\to\bc$ be the linear map with
$h_u(A_{\la})=P_{\la}(u)$ for each~$\la\in P^+$.

\begin{prop}\label{inc}(cf. \cite[Theorem 3.3.12]{macsph})
\begin{enumerate}
\item[(i)] Each algebra homomorphism $h:\sca\to\bc$ is of the form
$h=h_u$ for some $u\in\Hom(P,\bc^{\times})$, and
\item[(ii)] $h_{u}=h_{u'}$ if and only if $u'=wu$ for some $w\in W_0$.
\end{enumerate}
\end{prop}

\begin{proof}
Since $\bc[P]$ is integral over $\bc[P]^{W_0}$ \cite[V,
Exercise~12]{atiyah}, every homomorphism $\bc[P]^{W_0}\to\bc$ is
the restriction of a homomorphism $\bc[P]\to\bc$, and (i) follows.

It is clear that if $u'=wu$ for some $w\in W_0$, then
$h_{u'}=h_u$. Suppose now that $h_u=h_{u'}$. Since
$\{P_{\la}(x)\}_{\la\in P^+}$ forms a basis of $\bc[P]^{W_0}$ we see
that $u$ and $u'$ agree on $\bc[P]^{W_0}$, and thus their kernels
are maximal ideals of $\bc[P]$ lying over the same maximal ideal
of $\bc[P]^{W_0}$. The result now follows by \cite[V, Exercise
13]{atiyah}.
\end{proof}

We call the formula $h_u(A_{\la})=P_{\la}(u)$ the
\textit{Macdonald formula for the algebra homomorphisms
$\sca\to\bc$}. By comparing Proposition~\ref{newlem} and
\cite[Corollary~3.2.5]{macsph} we see that in the case when $P=Q$
(that is, when $R$ is of type $E_8$, $F_4$, $G_2$ or $BC_n$ for
some $n\geq1$) our formula $P_{\la}(u)$ agrees with the formula
in \cite[Theorem~4.1.2]{macsph}. The reason we require $P=Q$ here
is because in \cite{macsph} $u\in\Hom(Q,\bc^{\times})$ (although
$u$ there is called $s$), which is a consequence of Macdonald
requiring his group $G$ of $p$-adic type to be simply connected.

\begin{rem}\label{par} The function $P_{\la}(u)$ may be regarded as a
function of the variables $u_i=u^{\la_i}\in\bc^{\times}$,
$i=1,\ldots,n$. However, in general the coroots $\alpha^{\vee}$,
$\alpha\in R^+$, appearing in the formula for $c(u)$ do not have
particularly neat expressions in terms of the basis
$\{\la_i\}_{i=1}^n$. Thus in any given specific case it is often
useful to work with numbers other than the $\{u_i\}_{i=1}^{n}$.

Let us illustrate this in the $R=D_n$ case, which may be
described as follows (see \cite[Plate~IV]{bourbaki}). Let
$E=\br^n$ with standard orthonormal basis $\{e_i\}_{i=1}^{n}$, and
take $R$ to be the set of vectors $\pm e_i\pm e_j$, $1\leq
i<j\leq n$, where the $\pm$ signs may be taken independently. We
have $R^{\vee}=R$, and the set $\{e_i-e_{i+1},e_{n-1}+e_n\}_{1\leq
i\leq n-1}$ forms a base of $R$. The corresponding set of positive
roots is $\{e_i-e_j,e_i+e_j\}_{1\leq i<j\leq n}$. Observe that
$e_1=\la_1$, $e_i=\la_i-\la_{i-1}$ for $2\leq i\leq n-2$,
$e_{n-1}=\la_{n-1}+\la_n-\la_{n-2}$ and $e_n=\la_n-\la_{n-1}$.
Thus, defining numbers $t_i\in\bc^{\times}$, $i=1,\ldots,n$, by
$t_1=u_1$, $t_i=u_iu_{i-1}^{-1}$ $(2\leq i\leq n-2)$,
$t_{n-1}=u_{n-1}u_nu_{n-2}^{-1}$ and $t_n=u_nu_{n-1}^{-1}$, we
have
$$c(u)=\prod_{1\leq i<j\leq
n}\frac{(1-q^{-1}t_i^{-1}t_j)(1-q^{-1}t_i^{-1}t_j^{-1})}{(1-t_i^{-1}t_j)(1-t_i^{-1}t_j^{-1})}.$$
(Notice that $q_i=q$ for all $i\in I$, for there is only one root length).
\end{rem}

\section{The Integral Formula for the Algebra Homomorphisms}

We begin with some background that involves only the root system
$R$, which throughout is assumed to be irreducible. Define a
partial order on $P$ by setting $\mu\preceq\la$ if $\la-\mu\in
Q^+$.

We say that a subset $X\subset P$ is \textit{saturated} \cite[VI,
\ts1, Exercise~23]{bourbaki} if $\mu-i\alpha^{\vee}\in X$ for all
$\mu\in X$, $\alpha\in R$, and all $i$ between $0$ and
$\lan\mu,\alpha\ran$ inclusive. Every saturated set is stable
under $W_0$, and for each $\la\in P^+$ there is a unique saturated
set, denoted~$\Pi_{\la}$, with highest coweight~$\la$ (that is,
$\mu\preceq\la$ for all $\mu\in \Pi_{\la}$). We have
\begin{align}\label{eq:sat}
\Pi_{\la}=\{w\mu\mid\mu\in P^+, \mu\preceq\la,w\in W_0\}
\end{align}
(see \cite[Lemma~13.4B]{h2} for example).

We recall the definition of the \textit{Bruhat order} on $W$
\cite[\ts5.9]{h}. Let $v,w\in W$, and write $v\to w$ if
$v=s_{\alpha;k}w$ for some $\alpha\in R$, $k\in\bz$, and
$\ell(v)<\ell(w)$. Declare $v\leq w$ if and only if there exists
a sequence $v=w_0\to w_1\to\cdots\to w_n=w$. This gives the
Bruhat (partial) order on $W$. We extend the Bruhat order to
$\tilde{W}$ as in \cite[\ts2.3]{m} by declaring $\tilde{v}\leq
\tilde{w}$ if and only if $\tilde{v}=vg$ and $\tilde{w}=wg$ with
$v\leq w$ in $W$ and $g\in G$.

By a \textit{sub-expression} of a fixed reduced expression
$s_{i_1}\cdots s_{i_r}\in W$ we mean a product of the form
$s_{i_{k_1}}\cdots s_{i_{k_q}}$ where $1\leq k_1<\cdots<k_q\leq
r$. Let $w=s_{i_1}\cdots s_{i_r}$ be a fixed reduced expression
for $w\in W$. By \cite[Theorem~5.10]{h}, $v\leq w$ if and only $v$
can be obtained as a sub-expression of this reduced expression.

\begin{prop}\label{paper2prop2} Let $\tilde{v},\tilde{w}\in\tilde{W}$ with $\tilde{v}\leq \tilde{w}$.
If $\tilde{w}(0)\in \Pi_{\la}$, then $\tilde{v}(0)\in \Pi_{\la}$
too.
\end{prop}

\begin{proof} Suppose first that $\tilde{v}=s_{\alpha;k}\tilde{w}$
with $\ell(\tilde{v})<\ell(\tilde{w})$. Then by \cite[(2.3.3)]{m}
$H_{\alpha;k}$ separates $C_0$ and $\tilde{w}C_0$, and thus
$\lan\tilde{w}(0),\alpha\ran-k$ is between $0$ and
$\lan\tilde{w}(0),\alpha\ran$ (inclusive). Thus by the definition
of saturated sets
$$\tilde{v}(0)=\tilde{w}(0)-(\lan\tilde{w}(0),\alpha\ran-k)\alpha^{\vee}\in\Pi_{\la},$$
and the result clearly follows by induction.
\end{proof}

Let $\scx$ be an irreducible regular affine building. A
\textit{sector} of $\scx$ is a subcomplex $\cs\subset\scx$ such
that there exists an apartment $\ca$ such that $\cs\subset\ca$
and a type preserving isomorphism $\psi:\ca\to\S$ such that
$\psi(\cs)$ is a sector of $\S$. The \textit{base vertex} of
$\cs$ is $\psi^{-1}(\la)$, where $\la\in P$ is the base vertex of
$\psi(\cs)$.

If $\cs$ and $\cs'$ are sectors of $\scx$ with $\cs'\subseteq\cs$,
then we say that $\cs'$ is a \textit{subsector} of $\cs$. The
\textit{boundary} $\O$ of $\scx$ is the set of equivalence classes
of sectors, where we declare two sectors equivalent if and only if
they contain a common subsector. Given $x\in V_P$ and $\o\in\O$,
there exists a unique sector, denoted $\cs^x(\o)$, in the class
$\o$ with base vertex~$x$ \cite[Lemma~9.7]{ronan}.

\begin{lem}\label{isosect} Let $\cs$ be a sector in an apartment $\ca$ of $\scx$. There exists a unique type rotating isomorphism
$\psi_{\ca,\cs}:\ca\to\S$ such that $\psi_{\ca,\cs}(\cs)=\cs_0$.
\end{lem}

\begin{proof} Let $x$ be the base vertex of $\cs$, and let $\phi:\ca\to\S$ be
a type preserving isomorphism. Writing $\la=\phi(x)$ we see that
$t_{-\la}\circ\phi:\ca\to\S$ is a type rotating isomorphism
mapping $\cs$ to a sector of $\S$ based at $0$. Thus
$(t_{-\la}\circ\phi)(\cs)=w\cs_0$ for some $w\in W_0$, and so
$w^{-1}\circ t_{-\la}\circ\phi:\ca\to\S$ is a type rotating
isomorphism satisfying the requirements of the lemma.

Let $\psi$ and $\psi'$ be two such isomorphisms. It follows from
\cite[Lemma~5.3]{p} that $\psi'\circ\psi^{-1}=w$ for some $w\in
W_0$ (it is important that $\psi$ and $\psi'$ are type rotating
here). We have $w(C_0)=C_0$, and so $w=1$ since $W$ acts simply
transitively on the chambers of $\S$ \cite[p.142]{brown}. Thus
$\psi'=\psi$.
\end{proof}

Given an apartment $\ca$ and a sector $\cs$ of~$\ca$, let
$\rho_{\ca,\cs}:\scx\to\ca$ be the \textit{retraction onto $\ca$
centered at $\cs$} \cite[pages 170--171]{brown}. This is defined
as follows. Given any chamber $c$ of $\scx$, there exists a
subsector $\cs'$ of $\cs$ and an apartment $\ca'$ such that $c$
and $\cs'$ are contained in $\ca'$ \cite[page~170]{brown}. Writing
$\psi_{\ca'}:\ca'\to\ca$ for the (unique) isomorphism from
building axiom (iii), we set $\rho_{\ca,\cs}(c)=\psi_{\ca'}(c)$,
which is easily seen to be independent of the particular $\cs'$
and $\ca'$ chosen.

If $f=i_1\cdots i_r\in I^*$, the free monoid on $I$, we write
$s_f=s_{i_1}\cdots s_{i_r}\in W$. For $\la\in P^+$ let $w_{\la}$
be as in Section~\ref{uu}. Let $f_{\la}\in I^{*}$ be any (fixed)
word of minimal length such that $s_{f_{\la}}=w_{\la}$.

\begin{thm}\label{paper2thm1} Let $x\in V_P$, $\o\in \O$ and write $\cs=\cs^x(\o)$. Let $\ca$ be any apartment containing
$\cs$. Then $(\psi_{\ca,\cs}\circ\rho_{\ca,\cs})(y)\in\Pi_{\la}$
for all $y\in V_{\la}(x)$.
\end{thm}

\begin{proof} It follows easily from \cite[Proposition~5.2]{p} that there
exists a gallery $c_0,\ldots,c_n$ of type~$\s_i(f_{\la})$ from $x$
to $y$, where~$i=\tau(x)$. Write
$\Phi=\psi_{\ca,\cs}\circ\rho_{\ca,\cs}$. Then
$\Phi(c_0),\ldots,\Phi(c_n)$ is a pre-gallery of type $f_{\la}$
from $0$ to $\mu=\Phi(y)$. Thus there is a gallery of type
$f_{\la}'$, say, from $0$ to $\mu$, where $s_{f_{\la}'}$ is a
sub-expression of $s_{f_{\la}}$. Thus $s_{f_{\la}'}\leq
s_{f_{\la}}$, and so $s_{f_{\la}'}g_{l^{\vphantom{-1}}}\leq
s_{f_{\la}}g_{l^{\vphantom{-1}}}$, where $l=\tau(\la)$. Since
$(s_{f_{\la}}g_{l^{\vphantom{-1}}})(0)=\la\in\Pi_{\la}$, it
follows from Proposition~\ref{paper2prop2} that
$\mu'=(s_{f_{\la}'}g_{l^{\vphantom{-1}}})(0)\in\Pi_{\la}$ too.
Since $\Phi(c_0)=wC_0$ for some $w\in W_0$, we have
$\Phi(c_n)=ws_{f_{\la}'}C_0$, and so by considering types of
vertices we have
$$\mu=ws_{f_{\la}'}(g_{l^{\vphantom{-1}}}(0))=w\mu'.$$
Thus $\mu\in\Pi_{\la}$ by (\ref{eq:sat}).
\end{proof}

For each $x\in V_P$, $\o\in\O$ and $\la\in P^+$, the intersection
$V_{\la}(x)\cap\cs^x(\o)$ contains a unique vertex, denoted
$v_{\la}^x(\o)\in V_P$.

The coweights $h(x,y;\o)$ in the next theorem are the analogs of
the well studied \textit{horocycle numbers} of homogeneous trees.

\begin{thm}\label{c} Let $\o\in\O$ and let $x,y\in V_P$.
\begin{enumerate}
\item[(i)] Let $z\in\cs^x(\o)\cap\cs^y(\o)$ and write
$z=v_{\nu}^x(\o)=v_{\eta}^y(\o)$. The coweight $\nu-\eta$ is
independent of the particular $z\in\cs^x(\o)\cap\cs^y(\o)$
chosen. We denote this common value by $h(x,y;\o)$. If $\mu\in
P^+$ and $\mu-\nu\in P^+$ then
$$v_{\mu}^x(\o)=v_{\mu-h(x,y;\o)}^y(\o).$$
\item[(ii)] Suppose that $y\in V_{\la}(x)$. Write $\cs=\cs^x(\o)$ and let $\ca$ be any apartment containing $\cs$. Then
$h(x,y;\o)=(\psi_{\ca,\cs}\circ\rho_{\ca,\cs})(y)\in\Pi_{\la}.$
\end{enumerate}
\end{thm}

\begin{proof} (i) We have $v_{\nu+\mu'}^x(\o)=v_{\eta+\mu'}^y(\o)$ for all $\mu'\in
P^+$, since both are equal to $v_{\mu'}^z(\o)$. Thus, writing
$\mu=\mu'+\nu$ we have
\begin{align}
    \label{.1}
        v_{\mu}^x(\o)=v_{\mu-(\nu-\eta)}^y(\o)\qquad\textrm{whenever }\mu-\nu\in P^+.
\end{align}
If we instead choose $z'\in\cs^x(\o)\cap\cs^y(\o)$, where
$z'=v_{\nu'}^x(\o)=v_{\eta'}^y(\o)$, then following the above we
have
\begin{align}
    \label{.2}
        v_{\mu}^x(\o)=v_{\mu-(\nu'-\eta')}^y(\o)\qquad\textrm{whenever }\mu-\nu'\in P^+.
\end{align}
By choosing $\mu\in P^+$ such that both $\mu-\nu$ and $ \mu-\nu'$
are dominant, it follows from (\ref{.1}) and (\ref{.2}) that
$\nu-\eta=\nu'-\eta'$. Then (\ref{.1}) proves the final claim.

(ii) Write $\mu=(\psi_{\ca,\cs}\circ\rho_{\ca,\cs})(y)$. By
\cite[page~170]{brown} there exists a subsector $\cs'$ of
$\cs=\cs^x(\o)$ such that $\cs'$ and $y$ lie in a common
apartment $\ca'$. The restriction of $\rho_{\ca,\cs}$ to $\ca'$ is
thus a type preserving isomorphism. Pick $\nu\in P^+$ such that
$v_{\nu}^x(\o)\in\cs'$, $\nu-\mu\in P^+$ and
$v_{\nu}^x(\o)\in\cs^y(\o)$. The map
$\psi=t_{-\mu}\circ\psi_{\ca,\cs}\circ\rho_{\ca,\cs}:\ca'\to\S$ is
a type rotating isomorphism such that $\psi(y)=0$ and
$\psi(v_{\nu}^x(\o))=\nu-\mu\in P^+$. Thus $v_{\nu}^x(\o)\in
V_{\nu-\mu}(y)\cap \cs^y(\o)$, and so $h(x,y;\o)=\mu$. The fact
that $h(x,y;\o)\in\Pi_{\la}$ now follows from
Theorem~\ref{paper2thm1}.
\end{proof}

\begin{prop}\label{cocycle} For all $x,y,z\in V_P$ and $\o\in\O$ we
have the `\textit{cocycle relation}'
$h(x,y;\o)=h(x,z;\o)+h(z,y;\o)$. Thus $h(x,x;\o)=0$ and
$h(y,x;\o)=-h(x,y;\o)$.
\end{prop}

\begin{proof} For $\mu=k_1\la_1+\cdots+k_n\la_n\in P^+$ with each $k_i$ sufficiently large we have
$$v^x_{\mu}(\o)=v^z_{\mu-h(x,z;\o)}(\o)=v_{\mu-h(x,z;\o)-h(z,y;\o)}^y(\o)$$
and the result follows.
\end{proof}

The following theorem shows that if $y\in V_{\la}(x)$, then for
any $\o\in\O$, $\cs^y(\o)$ contains all vertices $v_{\mu}^x(\o)$
for $\mu\in P^+$ large enough, where \textit{large enough}
depends only on $\la$, not on the particular $x,y$ and $\o$.

For $\la\in P^+$, write $\mu\gg\la$ to mean that
$\mu-\Pi_{\la}\subset P^+$ (in particular, $\mu\in P^+$).

\begin{thm}\label{b} Let $x\in V_P$, $\la\in P^+$ and $y\in
V_{\la}(x)$. Then $v_{\mu}^x(\o)\in\cs^y(\o)$ for all $\o\in\O$
and all $\mu\gg\la$, and so
$v_{\mu}^x(\o)=v_{\mu-h(x,y;\o)}^y(\o)$ for all $\o\in\O$ and all
$\mu\gg\la$.
\end{thm}

\begin{rem} The proof of Theorem~\ref{b} will be given after the
following preliminary results. We are thankful to an anonymous
referee for sketching this proof of the present form of
Theorem~\ref{b}, which replaces our less sharp version of this
result.
\end{rem}

\begin{lem}\label{lem:restrictrho} Suppose that $\ca_1$ and~$\ca_2$ are apartments
containing a common sector~$\cs$. Then the maps
$\rho_{\ca_1,\cs}|_{\ca_2}:\ca_2\to \ca_1$ and
$\rho_{\ca_2,\cs}|_{\ca_1}:\ca_1\to \ca_2$ are mutually inverse
isomorphisms which fix~$\ca_1\cap \ca_2$ pointwise.
\end{lem}

\begin{proof} Fix any chamber $c\subset \cs$, and let
$\varphi:\ca_1\to\ca_2$ be the unique isomorphism fixing
$\ca_1\cap \ca_2$ pointwise. Then by definition we have
$\rho_{\ca_2,\cs}|_{\ca_1}=\varphi$, and since
$\varphi^{-1}:\ca_2\to\ca_1$ is the unique isomorphism fixing
$\ca_1\cap \ca_2$ pointwise we have
$\rho_{\ca_1,\cs}|_{\ca_2}=\varphi^{-1}$.
\end{proof}

\begin{lem}\label{lem:subsector} Let $\ca$ be an apartment in~$\scx$, let $\cs$ be
a sector in~$\ca$, and let $H$ be a wall in~$\ca$. Then exactly
one of the two closed half-apartments in~$\ca$ determined by~$H$,
$H^+$ say, contains a subsector of~$\cs$. If $x\in V_P\cap H^+$,
then the sector based at~$x$ and equivalent to~$\cs$ is contained
in~$H^+$.
\end{lem}
\begin{proof} Let $\psi=\psi_{\ca ,\cs}$ (see Lemma~\ref{isosect}). Then
$\psi(H)=H_{\alpha;k}$ for some $\alpha\in R$ and $k\in\bz$.
Since $H_{-\alpha;k}=H_{\alpha;-k}$, we may suppose that
$\alpha\in R^+$. If $k\leq 0$, then $\cs_0\subset H_{\alpha;k}^+=
\{x:\langle x,\alpha\rangle\geq k\}$. If $k\geq 1$, then
$\lambda+\cs_0\subset H_{\alpha;k}^+$ for
$\lambda=k(\lambda_1+\cdots+\lambda_n)$.

The final statement follows from \cite[Lemma~9.1]{ronan}.
\end{proof}

\goodbreak\begin{prop}\label{prop:nice} Let $c_0, ...,c_m$ be a gallery of type
$j_1\cdots j_m$, and let $\omega\in\Omega$. Let $\ca$ be an
apartment containing~$c_0$ and a sector~$\cs$ in the
class~$\omega$. Let $\rho=\rho_{\cs,\ca}$, and let
$e_k=\rho(c_k)$ for $k=0,\ldots,m$. For $k=1,\ldots,m$, let $H_k$
denote the wall in~$\ca$ containing the panel in~$e_{k-1}$ and
in~$e_k$ of type~$j_k$. Let $H_k^+$ denote the half-apartment
in~$\ca$ bounded by~$H_k$ which contains a subsector of~$\cs$
(see Lemma~\ref{lem:subsector}).

Then there exists an apartment $\cb$ containing $c_m$ and
$$\bigcap_{k=1}^m H_k^+$$
(and therefore $\cb$ contains a sector in the class~$\omega$).
\end{prop}
\begin{proof} By induction on~$m$. If $m=1$, the
panel of cotype~$j_1$ common to~$e_0=c_0$ and~$e_1$ is contained
in~$c_1$ and in~$H_1$. So by the proof of
\cite[Lemma~9.4]{ronan}, there is an apartment~$\cb$ containing
$H_1^+$ and~$c_1$.

Now suppose that $m>1$ and that there is an apartment $\cb'$
containing $c_{m-1}$ and $\bigcap_{i=1}^{m-1}H_k^+$. Let $\cs'$ be
any sector in the class~$\omega$ contained in
$\bigcap_{k=1}^{m-1}H_k^+$, and let $\rho'=\rho_{\cs',\cb'}$. If
$\ct$ is a common subsector of~$\cs$ and~$\cs'$, then
$\rho=\rho_{\ct,\ca}$ and $\rho'=\rho_{\ct,\cb'}$. So by
Lemma~\ref{lem:restrictrho}, the maps $\rho|_{\cb'}:\cb'\to \ca$
and $\rho'|_{\ca}:\ca\to \cb'$ are mutually inverse isomorphisms.

Let $H$ denote the wall in~$\cb'$ containing the panel of
cotype~$j_m$ in~$c_{m-1}$ and in~$c_m$, and let $H^+$ denote the
half-apartment in~$\cb'$ bounded by~$H$ and containing a subsector
of~$\cs$. The half-apartment $\rho'(H_m^+)$ in $\cb'$ is bounded
by~$H$. To see this, let $\pi$ denote the panel in $c_{m-1}$
and~$c_m$ of cotype~$j_m$, and let $\pi'$ denote the panel in
$e_{m-1}=\rho(c_{m-1})$ and~$e_m=\rho(c_m)$ of cotype~$j_m$. Then
$\pi'=\rho(c_{m-1}\cap c_m)$, and so $\rho'(\pi')=c_{m-1}\cap
c_m=\pi$. Now $\pi'$ is in the wall of~$\ca$ bounding~$H_m^+$,
and so $\pi=\rho'(\pi')$ is in the wall of~$\cb$ bounding
$\rho'(H_m^+)$. But $\pi$ is in the wall~$H$ of~$\cb$ bounding
$H^+$. Furthermore, $\rho'(H_m^+)$ contains a subsector of~$\cs$,
because $H_m^+\cap \cb'$ contains such a subsector, and is fixed
by~$\rho'$. Thus
\begin{equation}\label{eq:halfapts}
\rho'(H_m^+)=H^+.
\end{equation}

By the proof of \cite[Lemma~9.4]{ronan}, there is an apartment
$\cb$ containing~$c_m$ and~$H^+$. Since $H_1^+\cap\cdots\cap
H_{m-1}^+$ is contained in~$\cb'$, it is fixed by~$\rho'$, and so
$$H_1^+\cap\cdots\cap H_m^+ =\rho'((H_1^+\cap\cdots\cap
H_{m-1}^+)\cap H_m^+)\subset \rho'(H_m^+).$$ Thus by
(\ref{eq:halfapts}), $H_1^+\cap\cdots\cap H_m^+\subset H^+\subset
\cb$, completing the induction step.
\end{proof}

\begin{lem}\label{lem:bruhat} Suppose that $C_0=D_0,\ldots,D_m$ is a gallery in $\Sigma$ joining $0$ to~$\la$ of minimal
length, and that $C_0=E_0,\ldots,E_m$ is a pre-gallery
in~$\Sigma$. If both the gallery and the pre-gallery have the
same type, then each type $l=\tau(\la)$ vertex of any of the
$D_i$'s and~$E_i$'s is in~$\Pi_\lambda$.
\end{lem}
\begin{proof} For $k=0,\ldots,m$, define $u_k,v_k\in W$ by
$D_k=u_kC_0$ and $E_k=v_kC_0$. Then the type $l$ vertices of
$D_k$ and $E_k$ are $u_kg_l(0)$ and $v_kg_l(0)$ respectively (see
the proof of Theorem~\ref{paper2thm1}). In particular,
$u_mg_l(0)=\la\in\Pi_{\la}$.

Since the gallery $C_0=D_0,\ldots,D_m$ and the pre-gallery
$C_0=E_0,\ldots,E_m$ have the same type, $j_1\ldots j_m$, say,
each $u_k$ and each~$v_k$ is a subexpression of the reduced
expression $u_m=s_{j_1}\cdots s_{j_m}$. Thus $u_kg_l,v_kg_l\leq
u_mg_l$ (with respect to the Bruhat order on $\tilde{W}$) for each
$k=0,\ldots,m$.

The result now follows from Proposition~\ref{paper2prop2}.
\end{proof}

\begin{cor}\label{cor:nice} Suppose that $x,y\in V_P$, $\o\in\O$, $\la\in P^+$, and $y\in V_{\la}(x)$. Consider a minimal gallery
$c_0,...,c_m$ from $x$ to~$y$ and an apartment~$\ca$ containing
$c_0$ and a sector in the class~$\omega$ (and hence
containing~$\cs^x(\omega)$). Let $\psi:\ca\to\Sigma$ be the
type-rotating isomorphism mapping $\cs^x(\omega)$ to~$\cs_0$ (see
Lemma~\ref{isosect}). Finally, let $H_i$ and $H_i^+$ be as in
Proposition~\ref{prop:nice}, and write
$\rho=\rho_{\ca,\cs^x(\omega)}$. Then:
\begin{enumerate}
\item[(i)] $\psi^{-1}(\Pi_\lambda)$ contains all the type $j=\tau(y)$
vertices in each $\rho(c_i)$, and so in particular it contains
$y'=\rho(y)$.
\item[(ii)] If $\mu\gg\la$, then
$v^x_{\mu}(\omega)\in \bigcap_{i=1}^m  H_i^+$.
\end{enumerate}
\end{cor}

\begin{proof} (i) Let $e_i=\rho(c_i)$ and $E_i=\psi(e_i)$ for
$i=0,\ldots,m$. Let $\rho'=\rho_{c_0,\ca}$ denote the retraction
of center~$c_0$ onto~$\ca$ (see \cite[\ts IV.3]{brown}). Let
$d_i=\rho'(c_i)$ and $D_i=\psi(d_i)$ for $i=0,\ldots,m$. Then the
gallery $C_0=D_0,\ldots,D_m$ and the pre-gallery
$C_0=E_0,\ldots,E_m$ satisfy the hypotheses of
Lemma~\ref{lem:bruhat}, and the result follows.

(ii) For $i=0,\ldots,m$, let $v_i$ be the type~$j$ vertex
of~$e_i$, and so $\{v_i\}_{i=0}^m\subset\psi^{-1}(\Pi_{\la})$.

Let $1\leq i\leq m$. If $v_{i-1}=v_i$, then
$\psi(v_i)\in\Pi_\lambda\cap\psi(H_i)$, and so $\psi(v_i)+\cs_0$
is contained in~$\psi(H_i)^+$ by Lemma~\ref{lem:subsector} (here
$\psi(H_i)^+$ is the half-space of $\S$ bounded by $\psi(H_i)$
and containing a subsector of $\cs_0=\psi(\cs^x(\o))$). Hence
$\cs^{v_i}(\omega)$ is contained in~$H_i^+$. By our hypothesis,
$\mu-\Pi_\lambda\subset P^+$, and so $\mu=\psi(v_i)+\nu_i$ for
some $\nu_i\in P^+$. Hence
$\psi(v^x_\mu(\omega))=\mu\in\psi(v_i)+\cs_0\subset\psi(H_i)^+$.
Hence $v^x_\mu(\omega)\in H_i^+$.

If $v_{i-1}\neq v_i$, then $\psi(v_{i-1})$ and~$\psi(v_{i})$ lie
on opposite sides of~$\psi(H_i)$. Then by
Lemma~\ref{lem:subsector}, either $\psi(v_{i-1})+\cs_0$
or~$\psi(v_i)+\cs_0$ is contained in $\psi(H_i)^+$. Let us assume
that $\psi(v_i)+\cs_0\subset\psi(H_i)^+$. Since
$\mu-\Pi_\lambda\subset P^+$, we can write $\mu=\psi(v_i)+\nu_i$
for some $\nu_i\in P^+$. Hence
$\psi(v^x_\mu(\omega))=\mu\in\psi(v_i)+\cs_0\subset\psi(H_i)^+$,
and so again $v^x_\mu(\omega)\in H_i^+$.
\end{proof}

\begin{proof}[Proof of Theorem~\ref{b}] Let $c_0,\ldots,c_m$ be a
minimal gallery from~$x$ to~$y$, and let $\ca,\cb$ and
$H_1^+,\ldots,H_m^+$ be as in Proposition~\ref{prop:nice}.

By Lemma~\ref{lem:restrictrho}, the map
$\varphi=\rho_{\ca,\cs^x(\o)}|_{\cb}:\cb\to \ca$ is an
isomorphism. It maps~$y$ to~$y'=\rho_{\ca,\cs^x(\omega)}(y)$ and
fixes $\ca\cap \cb\supset\bigcap_{i=1}^mH_i^+$, which contains a
sector in the class~$\omega$. Hence $\varphi$ maps
$\cs^y(\omega)$ to~$\cs^{y'}(\omega)$. Moreover, if $\mu\gg\la$,
then by~Corollary~\ref{cor:nice}(ii), $\varphi$ fixes
$v^x_{\mu}(\omega)$.

Let $\psi:\ca\to\Sigma$ be the type-rotating isomorphism mapping
$\cs^x(\omega)$ to~$\cs_0$. Then by Corollary~\ref{cor:nice}(i) we
have $y'\in\psi^{-1}(\Pi_\lambda)$, and so we can write
$\mu=\psi(y')+\nu$ for some $\nu\in P^+$. Therefore
$v^x_\mu(\omega)\in\cs^{y'}(\omega)$, and applying~$\varphi$, we
see that $v^x_\mu(\omega)\in\cs^y(\omega)$. So
$v^x_\mu(\omega)=v^y_{\mu-h(x,y;\omega)}(\omega)$ by the
definition of~$h(x,y;\omega)$.
\end{proof}

Let $\leq$ denote the partial order on $P^+$ given by $\mu\leq
\la$ if and only if $\la-\mu\in P^+$. Fixing $x\in V_P$, there is
a natural map $\theta:\O\to\prod_{\la\in P^+}V_{\la}(x)$, where
one maps $\o$ to $(v_{\la}^x(\o))_{\la\in P^+}$. For each pair
$\la,\mu\in P^+$ with $\mu\leq\la$, let
$\varphi_{\mu,\la}:V_{\la}(x)\to V_{\mu}(x)$ be the map $y\mapsto
v_{\mu}(x,y)$, where $v_{\mu}(x,y)$ is the unique vertex in
$V_{\mu}(x)\cap\mathrm{conv}\{x,y\}$ (see Appendix~\ref{B}). Then
$(V_{\la}(x),\varphi_{\mu,\la})$ is an inverse system of
topological spaces (where each finite set $V_{\la}(x)$ is given
the discrete topology). The inverse limit
$\varprojlim(V_{\la}(x),\varphi_{\mu,\la})$ is a compact
Hausdorff topological space \cite[I.9.6]{bourbakisets}, and the
map $\theta$ is a bijection of $\O$ onto this inverse limit, thus
inducing a compact Hausdorff topology on $\O$, which we show in
Theorem~\ref{topology} is independent of $x\in V_P$. See
Appendix~\ref{B} for a sketch of the proof that $\theta$ is a
bijection of $\O$ onto
$\varprojlim(V_{\la}(x),\varphi_{\mu,\la})$.

With $x\in V_P$ fixed as above, for each $y\in V_P$ let
$\O_x(y)=\{\o\in\O\mid y\in\cs^x(\o)\}$. The sets $\O_x(y)$,
$y\in V_P$, form a basis of open and closed sets for the topology
on $\O$, and the functions $\o\mapsto h(x,y;\o)$ are locally
constant on $\O$, as we see in Lemma~\ref{d}.

To each $x\in V_P$ there is a unique regular Borel probability
measure $\nu_x$ on $\O$ such that $\nu_x(\O_x(y))=N_{\la}^{-1}$
if $y\in V_{\la}(x)$. To see this, for each $\la\in P^+$ let
$\scc_{\la}(\O)$ be the space of all functions $f:\O\to\bc$ which
are constant on each set $\O_x(y)$, $y\in V_{\la}(x)$. For each
$\la\in P^+$ define $J_{\la}:\scc_{\la}(\O)\to\bc$ by
$J_{\la}(f)=\frac{1}{N_{\la}}\sum_{y\in V_{\la}(x)}c_{y}(f)$,
where $c_{y}(f)$ is the constant value $f$ takes on $\O_x(y)$. The
space of all locally constant functions $f:\O\to\bc$ is
$\scc_{\infty}(\O)=\bigcup_{\la\in P^+}\scc_{\la}(\O)$. Define
$J:\scc_{\infty}(\O)\to\bc$ by $J(f)=J_{\la}(f)$ if $f\in
\scc_{\la}(\O)$. The map $J$ is linear, maps
$1\in\scc_{\infty}(\O)$ to $1\in\bc$ and satisfies
$|J(f)|\leq\|f\|_{\infty}$ for all $f\in \scc_{\infty}(\O)$. Since
$\scc_{\infty}(\O)$ is dense in $\scc(\O)$, $J$ extends uniquely
to a linear map $\tilde{J}:\scc(\O)\to\bc$ such that
$|\tilde{J}(f)|\leq \|f\|_{\infty}$ for all $f\in\scc(\O)$ (here
$\scc(\O)$ is the space of all continuous functions
$f:\O\to\bc$). Thus by the Riesz Representation Theorem there
exists a unique regular Borel probability measure $\nu_x$ such
that
$$\tilde{J}(f)=\int_{\O}f(w)d\nu_x(\o)\quad\textrm{for all
$f\in\scc(\O)$}.$$ In particular, $N_{\la}^{-1}=\nu_{x}(\O_x(y))$
if $y\in V_{\la}(x)$.

In Theorem~\ref{topology}(ii) we show that for $x,y\in V_P$, the
measures $\nu_x$ and $\nu_y$ are mutually absolutely continuous,
and we compute the \textit{Radon-Nikodym derivative}.

\begin{lem}\label{d} Let $y\in V_{\nu}(x)$ and suppose that $z\in
V_{\la}(x)\cap V_{\mu}(y)$ with $\la\gg\nu$. Then
\begin{enumerate}
\item[$(i)$] $\O_x(z)\subset\O_y(z)$, and
\item[$(ii)$] $h(x,y;\o)=\la-\mu$ for all $\o\in\O_x(z)$.
\end{enumerate}
\end{lem}

\begin{proof} (i) Let $\o\in\O_x(z)$, and so $z=v_{\la}^x(\o)$, and by Theorem~\ref{b} $z\in\cs^y(\o)$.
Thus $\o\in\O_y(z)$ and so
$\O_x(z)\subset\O_y(z)$. Note that if $\mu\gg\nu^*$ too, then
$\O_x(z)=\O_y(z)$.

Since $\O_x(z)\subset\O_y(z)$ and $z\in V_{\la}(x)\cap
V_{\mu}(y)$ we have $v_{\la}^x(\o)=z=v_{\mu}^y(\o)$ for all $\o\in
\O_x(z)$, and so (ii) follows from Theorem~\ref{c}(i).
\end{proof}

\begin{lem}\label{dd} Let $x,y\in V_P$. In the notation of Appendix~\ref{B}, if
$z\in\mathrm{conv}\{x,y\}$, then $\O_x(y)\subset \O_x(z)$.
\end{lem}

\begin{proof} Let $\o\in \O_x(y)$. Then the sector $\cs^x(\o)$
contains $x$ and $y$, and hence $z$. Thus $\o\in\O_x(z)$.
\end{proof}

Recall that we write $P^{++}$ for the set of all \textit{strongly
dominant coweights of $R$}, that is, those $\la\in P$ such that
$\lan\la,\alpha_i\ran>0$ for all $i\in I_0$.

\begin{prop}\label{topt} Let $\mu\in P$ be fixed. For all $\la\in P^{++}$
such that $\la-\mu\in P^{++}$ we have
$$\frac{N_{\la}}{N_{\la-\mu}}=\prod_{i=1}^{n}q_{t_{\la_i}}^{\lan\mu,\alpha_i\ran}=\prod_{\alpha\in
R^+}\tau_{\alpha}^{\lan \mu,\alpha\ran}$$
\end{prop}

\begin{proof} By (\ref{ns}) we have
$N_{\la}N_{\la-\mu}^{-1}=q_{t_{\la}}q_{t_{\la-\mu}}^{-1}$, and so
by (\ref{ns2})
$$\frac{N_{\la}}{N_{\la-\mu}}=\frac{q_{t_{\la}}}{q_{t_{\la-\mu}}}=\prod_{i=1}^n q_{t_{\la_i}}^{\lan\la,\alpha_i\ran-\lan\la-\mu,\alpha_i\ran}=\prod_{i=1}^{n}q_{t_{\la_i}}^{\lan\mu,\alpha_i\ran},$$
proving the first equality. On the other hand, by
Proposition~\ref{newlem},
\begin{align*}
\frac{N_{\la}}{N_{\la-\mu}}&=\frac{q_{t_{\la}}}{q_{t_{\la-\mu}}}=\prod_{\alpha\in
R^+}\tau_{\alpha}^{\lan\la,\alpha\ran-\lan\la-\mu,\alpha\ran}=\prod_{\alpha\in
R^+}\tau_{\alpha}^{\lan\mu,\alpha\ran}.\qedhere
\end{align*}
\end{proof}

\goodbreak
\begin{lem}\label{vs} Let $\la,\mu\in P^+$. Then
$\Pi_{\la}+\Pi_{\mu}\subseteq\Pi_{\la+\mu}$.
\end{lem}

\begin{proof} Let $\la'\in\Pi_{\la}$ and $\mu'\in\Pi_{\mu}$. Then
$$w(\la'+\mu')=w\la'+w\mu'\preceq\la+\mu\quad\textrm{for all
$w\in W_0$}.$$ By choosing $w\in W_0$ such that $w(\la'+\mu')\in
P^+$ we have $w(\la'+\mu')\in\Pi_{\la+\mu}$ by~(\ref{eq:sat}), and
so $\la'+\mu'\in\Pi_{\la+\mu}$.
\end{proof}

\begin{thm}\label{topology} Consider the topologies and measures defined above
on $\O$. Then
\begin{enumerate}
\item[(i)] The topology on $\O$ does not depend on the particular $x\in
V_P$.
\item[(ii)] For $x,y\in V_P$, the measures $\nu_x$ and
$\nu_y$ are mutually absolutely continuous, and the Radon-Nikodym
derivative is given by
$$\frac{d\nu_y}{d\nu_x}(\o)=\prod_{i=1}^{n}q_{t_{\la_i}}^{\lan
h(x,y;\o),\alpha_i\ran}=\prod_{\alpha\in R^+}\tau_{\alpha}^{\lan
h(x,y;\o),\alpha\ran}.$$
\end{enumerate}
\end{thm}

\begin{proof} (i) Let $x,y\in V_P$, with $y\in V_{\nu}(x)$, say, and choose
$\la\gg\nu+\nu^*$. Notice that this implies that $\la\gg\nu$.
Furthermore, for each $\nu'\in\Pi_{\nu}$, using Lemma~\ref{vs} we
have
$$(\la-\nu')-\Pi_{\nu^*}\subseteq\la-\Pi_{\nu}-\Pi_{\nu^*}\subseteq\la-\Pi_{\nu+\nu^*},$$
and so $\la-\nu'\gg\nu^*$.

Let $\o_0\in\O_x(v)$, a basic open set for the topology using $x$,
where $v\in V_{\eta}(x)$, say. Since $\la\gg\nu$, by
Theorem~\ref{b} we have $v_{\la}^x(\o_0)\in\cs^y(\o_0)$. Let
$z=v_{\la}^x(\o_0)=v_{\la-h(x,y;\o_0)}^y(\o_0)$, and so $z\in
V_{\la}(x)\cap V_{\la-h(x,y;\o_0)}(y)$. Now
$h(x,y;\o_o)\in\Pi_{\nu}$ (see Theorem~\ref{c}(ii)), and so by the
above, $\la-h(x,y;\o_0)\gg\nu^*$. Since $\la\gg\nu$ too, by
Lemma~\ref{d} we have $\O_x(z)=\O_y(z)$. Choosing $\la$ so that
$\la-\eta\in P^+$ we have $\O_x(z)\subset\O_x(v)$ by
Lemma~\ref{dd}. Thus, since $z=v_{\la-h(x,y;\o_0)}^y(\o_0)$, we
have $\o_0\in\O_y(z)=\O_x(z)\subset\O_x(v)$, which shows that the
$x$-\textit{open} sets are $y$-\textit{open}, and the first
statement follows.

(ii) With $x,y,z$ and $\o_0$ as above, since $\O_x(z)=\O_y(z)$,
and $\nu_x(\O_x(z))=N_{\la}^{-1}$, and
$\nu_y(\O_y(z))=N_{\la-h(x,y;\o_0)}^{-1}$, we see that
$$\nu_y(\O_x(z))=\nu_y(\O_y(z))=N_{\la-h(x,y;\o_0)}^{-1}=N_{\la-h(x,y;\o_0)}^{-1}N_{\la}\nu_x(\O_x(z)),$$
and so the measures are mutually absolutely continuous, and the
Radon-Nikodym derivative is given by
$$\frac{d\nu_y}{d\nu_x}(\o)=\frac{N_{\la}}{N_{\la-h(x,y;\o)}}\quad\textrm{for any $\la\in P^+$ such that
$\la\gg\nu+\nu^*$}.$$ The result follows from
Proposition~\ref{topt} by choosing $\la$ perhaps larger still so
that both $\la$ and $\la-h(x,y;\o)$ are strongly dominant.
\end{proof}

Let $r\in\Hom(P,\bc^{\times})$ be the map
\begin{align}
\label{r}\mu\mapsto\prod_{i=1}^{n}q_{t_{\la_i}}^{\frac{1}{2}\lan\mu,\alpha_i\ran}=\prod_{\alpha\in
R^+}\tau_{\alpha}^{\frac{1}{2}\lan \mu,\alpha\ran}
\end{align}
Following our usual convention we write $r^{\mu}$ in place of
$r(\mu)$.

Proposition~\ref{topt} immediately gives the following.

\begin{cor}\label{top} Let $\mu\in P$ be fixed. Then
$$r^{\mu}=\left(\frac{N_{\la}}{N_{\la-\mu}}\right)^{1/2}$$ for any $\la\in
P^{++}$ such that $\la-\mu\in P^{++}$.
\end{cor}

For $\la\in P^+$, let us write $\mu\ggg\la$ to mean that
$\mu-\Pi_{\la}\subset P^{++}$ (in particular, notice that if
$\mu\ggg\la$, then $\mu\gg\la$ and $\mu\in P^{++}$). The reason
for introducing this notation is that we will want to ensure that
the formula in Corollary~\ref{top} is applicable in the following
results.

Recall the definition of the numbers $a_{\la,\mu;\nu}$ from
before (\ref{a}).

\begin{lem}\label{paper2lem4} Let $\la\in P^+$. For each $x\in V_P$, $\o\in\O$, $\mu\in\Pi_{\la}$ and $\nu\ggg\la$,
$$\frac{1}{N_{\la}}|\{y\in V_{\la}(x)\mid
h(x,y;\o)=\mu\}|=r^{-2\mu}a_{\la,\nu-\mu;\nu}.$$ In particular,
the value of the left hand side is independent of $x\in V_{P}$
and $\o\in \O$.
\end{lem}

\begin{proof} We will first show that whenever $\nu\gg\la$,
\begin{align}\label{eq:mis}
\{y\in V_{\la}(x)\mid h(x,y;\o)=\mu\}=V_{\la}(x)\cap
V_{(\nu-\mu)^*}(v_{\nu}^x(\o)).
\end{align}
If $y\in V_{\la}(x)\cap V_{(\nu-\mu)^*}(v_{\nu}^x(\o))$, then by
Theorem~\ref{b}, $v_{\nu}^x(\o)\in\cs^y(\o)\cap V_{\nu-\mu}(y)$,
and so $v_{\nu}^x(\o)=v_{\nu-\mu}^y(\o)$. Thus $h(x,y;\o)=\mu$.

Conversely, if $y\in V_{\la}(x)$ and $h(x,y;\o)=\mu$, then
$v_{\nu}^x(\o)=v_{\nu-\mu}^y(\o)$ once $\nu\gg\la$ by
Theorem~\ref{b}. Thus $y\in V_{\la}(x)\cap
V_{(\nu-\mu)^*}(v_{\nu}^x(\o))$.

Now suppose that $\nu\ggg\la$. By (\ref{a}), (\ref{eq:mis}) and
Corollary~\ref{top} we have
\begin{align*}
\frac{1}{N_{\la}}|\{y\in V_{\la}(x)\mid
h(x,y;\o)=\mu\}|&=\frac{N_{\nu-\mu}}{N_{\nu}}a_{\la,\nu-\mu;\nu}=r^{-2\mu}a_{\la,\nu-\mu;\nu}.\qedhere
\end{align*}
\end{proof}

We now describe the algebra homomorphisms $h:\sca\to\bc$ in terms
of \textit{(zonal) spherical functions}.
\begin{defn}\label{zsf} Fix a vertex $x\in V_P$. A function
$f:V_P\to\bc$ is called \textit{spherical with respect to $x$} if
\begin{enumerate}
\item[(i)] $f(x)=1$,
\item[(ii)] $f$ is $x$-radial (that is, $f(y)=f(y')$ whenever $y,y'\in V_{\la}(x)$), and
\item[(iii)] for each $A\in\sca$ there is a number $c_A$ such that
$Af=c_Af$.
\end{enumerate}
\end{defn}

The following is proved in \cite[Proposition~3.4]{cm} in the
$\widetilde{A}_2$ case, and the proof there generalises
immediately.

\begin{prop}\label{conn} An $x$-radial function $f:V_P\to\bc$ is spherical if
and only if the map $h:\sca\to\bc$ given by $h(A)=(Af)(x)$
defines an algebra homomorphism. Moreover, each
$h\in\Hom(\sca,\bc)$ arises in this way.
\end{prop}

Let $x\in V_P$ be fixed and let $u\in\Hom(P,\bc^{\times})$ and
$y\in V_P$. Define
$$F^{x}_{u}(y)=\int_{\O}(ur)^{h(x,y;\o)}d\nu_x(\o),$$
where $(ur)^{\la}=u^{\la}r^{\la}$ for all $\la\in P$. The integral
exists by Theorem~\ref{c}(ii) and the fact that $N_{\la}$ and
$|\Pi_{\la}|$ are finite for each $\la\in P^+$.

In the following theorem we provide a second formula $h'_u$,
$u\in\Hom(P,\bc^{\times})$, for the algebra homomorphisms
$\sca\to\bc$. In Theorem~\ref{mainresult} we will show that
$h_u'=h_u$.

\begin{thm}\label{main} Let $x,y\in V_P$ with $y\in V_{\la}(x)$.
Then for all $u\in\Hom(P,\bc^{\times})$
\begin{itemize}
\item[(i)] $F^{x}_{u}(x)=1$,
\item[(ii)] $F^{x}_{u}(y)=F^{x'}_{u}(y')$ whenever $x',y'\in V_P$ satisfy $y'\in
V_{\la}(x')$, and
\item[(iii)] $A_{\la}F^{x}_{u}=\varphi_{\la}(u)F^x_{u}$, where
for any $\nu\ggg\la$
$$\varphi_{\la}(u)=\sum_{\mu\in\Pi_{\la}}r^{-\mu}a_{\la,\nu-\mu;\nu}u^{\mu},$$
which is independent of $x\in V_P$.
\end{itemize}
Thus the map $h_u':\sca\to\bc$ given by $h_u'(A)=(AF^x_u)(x)$
defines an algebra homomorphism (by Proposition~\ref{conn}).
\end{thm}

\begin{proof} Since $\nu_x$ is a probability measure, (i) follows from Proposition~\ref{cocycle}.

We now prove (ii), which we note is stronger than the claim that
$F^x_u$ is $x$-radial. Let $\nu\in P^+$. Since $\O$ is the union
of the disjoint sets $\O_x(z)$ over $z\in V_{\nu}(x)$, we have
\begin{align*}
F_{u}^{x}(y)&=\sum_{z\in V_{\nu}(x)}\int_{\O_x(z)}(ur)^{h(x,y;\o)}d\nu_x(\o)\\
&=\sum_{\mu\in P^+}\sum_{z\in V_{\nu}(x)\cap
V_{\mu}(y)}\int_{\O_x(z)}(ur)^{h(x,y;\o)}d\nu_x(\o).
\end{align*}
Now take $\nu\gg\la$, and so by Lemma~\ref{d} $h(x,y;\o)=\nu-\mu$
for all $\o\in\O_x(z)$ and $z\in V_{\nu}(x)\cap V_{\mu}(y)$. Since
$\nu_x(\O_x(z))=N_{\nu}^{-1}$ we have
$$F_{u}^x(y)=\sum_{\mu\in P^+}\frac{1}{N_{\nu}}|V_{\nu}(x)\cap
V_{\mu}(y)|(ur)^{\nu-\mu},$$ and the result follows from (\ref{a}).

We now prove (iii). Let $\nu\ggg\la$. By the cocycle relations
(Proposition~\ref{cocycle}), Theorem~\ref{c}(ii) and
Lemma~\ref{paper2lem4} we have
\begin{align*}
(A_{\la}F_{u}^x)(y)&=\frac{1}{N_{\la}}\sum_{z\in
V_{\la}(y)}\int_{\O}(ur)^{h(x,z;\o)}d\nu_x(\o)\\
&=\int_{\O}\bigg(\frac{1}{N_{\la}}\sum_{z\in
V_{\la}(y)}(ur)^{h(y,z;\o)}\bigg)(ur)^{h(x,y;\o)}d\nu_x(\o)\\
&=\bigg(\sum_{\mu\in\Pi_{\la}}r^{-\mu}a_{\la,\nu-\mu;\nu}u^{\mu}\bigg)\int_{\O}(ur)^{h(x,y;\o)}d\nu_x(\o).\qedhere
\end{align*}
\end{proof}

\begin{cor}\label{2ndform} Let $y\in V_{\la}(x)$. Then for any $\o\in\O$,
$$h'_u(A_{\la})=F_{u}^x(y)=\frac{1}{N_{\la}}\sum_{z\in
V_{\la}(x)}(ur)^{h(x,z;\o)}=\varphi_{\la}(u).$$
\end{cor}

\begin{proof} By Theorem~\ref{main}(ii) and the definition of $A_{\la}$ we have
$(A_{\la}F^x_{u})(x)=F^x_{u}(y)$, and by Theorem~\ref{main}(i) we
have $\varphi_{\la}(u)F^x_{u}(x)=\varphi_{\la}(u)$. The result
now follows from Theorem~\ref{main}(iii) and the proof thereof.
\end{proof}

\section{The Plancherel Measure} Each $A\in\sca$ maps $\ell^2(V_P)$ into itself, and for $\la\in P^+$ and
$f\in\ell^2(V_P)$ we have $\|A_{\la}f\|_2\leq\|f\|_2$ (see
\cite[Lemma~4.1]{C2} for a proof in a similar context). So we may
regard $\sca$ as a subalgebra of the $C^*$-algebra
$\mathscr{L}(\ell^2(V_P))$ of bounded linear operators on
$\ell^2(V_P)$. The facts that $y\in V_{\la}(x)$ if and only if
$x\in V_{\la^*}(y)$, and $N_{\la^*}=N_{\la}$, imply that
$A_{\la}^*=A_{\la^*}$, and so the adjoint $A^*$ of any $A\in\sca$
is also in~$\sca$.

Let $\sca_2$ denote the completion of $\sca$ with respect to
$\|\cdot\|$, the $\ell^2$-operator norm. So $\sca_2$ is a
commutative $C^*$-algebra. We write $M_2=\Hom(\sca_2,\bc)$ (this
is the \textit{maximal ideal space of $\sca_2$}), and we denote
the associated Gelfand transform $\sca_2\to\scc(M_2)$ by
$A\mapsto\hat{A}$, where $\hat{A}(h)=h(A)$. Here $\scc(M_2)$ is
the algebra of $\mathrm{w}^*$-continuous functions on $M_2$ with
the sup norm. This map is an isometric isomorphism of
$C^*$-algebras \cite[Theorem~I.3.1]{davidson}.

The algebra homomorphisms $\tilde{h}:\sca_2\to\bc$ are precisely
the extensions to $\sca_2$ of algebra homomorphisms $h:\sca\to\bc$
which are continuous with respect to the $\ell^2$-operator norm.
When there is no risk of ambiguity we will simply write $h$ in
place of $\tilde{h}$. If $h=h_u:\sca_2\to\bc$ we write
$\hat{A}(u)$ in place of $\hat{A}(h)$ (so
$\widehat{A}_{\la}(u)=P_{\la}(u)$).

Let $\lan\cdot,\cdot\ran$ be the usual inner product on
$\ell^2(V_P)$ (this is not to be confused with the unrelated inner
product on $E$). If $X\subset V_P$ we write $1_X$ for the
characteristic function on $X$, and we write $\delta_x$ for
$1_{\{x\}}$.
\begin{lem}\label{inner} Let $A\in\sca_2$ and $o\in V_P$. Then $A\delta_o=0$
implies that $A=0$.
\end{lem}

\begin{proof} Let $x\in V_P$. Observe that if $A\in\sca$ then $A\delta_x$ is
$x$-radial, for if $A=\sum_{\la}a_{\la}A_{\la}$ is a finite linear
combination, then
$A\delta_x=\sum_{\la}a_{\la}N_{\la}^{-1}1_{V_{\la^*}(x)}$, which
is $x$-radial. It follows that $A\delta_x$ is $x$-radial for all
$A\in\sca_2$. Now, given $A\in\sca_2$ and $\mu\in P^+$, $\lan
A\delta_x,1_{V_{\mu^*}(x)}\ran$ does not depend on $x\in V_P$, for
if $A=\sum_{\la}a_{\la}A_{\la}\in\sca$, then $\lan
A\delta_x,1_{V_{\mu^*}(x)}\ran=a_{\mu}$. Thus if $A\in\sca_2$ and
$A\delta_o=0$, then $\lan A\delta_o,1_{V_{\mu^*}(o)}\ran=0$ for
all $\mu\in P^+$, and so $\lan A\delta_x,1_{V_{\mu^*}(x)}\ran=0$
for all $\mu\in P^+$ and for all $x\in V_P$. Since $A\delta_x$ is
$x$-radial for all $x\in V_P$, it follows that $A\delta_x=0$ for
all $x$, and so $Af=0$ for all finitely supported functions
$f\in\ell^2(V_P)$. Thus by density the same is true for all $f\in
\ell^2(V_P)$, completing the proof.
\end{proof}

Since $A_{\la}\delta_o=N_{\la}^{-1}1_{V_{\la^*}(o)}$ for each
$o\in V_P$, we have $\lan
A_{\la}\delta_o,A_{\mu}\delta_o\ran=\delta_{\la,\mu}N_{\la}^{-1}$.
Thus by Lemma~\ref{inner} we can define an inner product on
$\sca_2$ (independent of $o\in V_P$) by $\lan A,B\ran=\lan
A\delta_o,B\delta_o\ran$.

For any fixed $o\in V_P$, the map $A\mapsto\lan
A,A_0\ran=(A\delta_o)(o)$ maps the identity $A_0$ of $\sca_2$ to 1
and satisfies $|\lan
A,A_0\ran|\leq\|A\delta_o\|_2\leq\|A\|=\|\hat{A}\|_{\infty}$.
Thus by the Riesz Representation Theorem there exists a unique
regular Borel probability measure $\pi$ on $M_2$ so that
$$(A\delta_o)(o)=\int_{M_2}\hat{A}(h)d\pi(h)\quad\textrm{for all
$A\in\sca_2$.}$$ Hence, for all $A,B\in\sca_2$,
\begin{align}
\label{plancherelmeasure}\lan A,B\ran=\lan A
\delta_o,B\delta_o\ran=(B^*A\delta_o)(o)=\int_{M_2}\hat{A}(h)\overline{\hat{B}(h)}d\pi(h).
\end{align}
We refer to $\pi$ and $M_2$ as the \textit{Plancherel measure}
and \textit{spectrum} of $\sca_2$, respectively.

\begin{prop}\label{plancherelthm} $M_2=\mathrm{supp}(\pi)$.
\end{prop}

\begin{proof} If $h_0\in M_2\backslash\mathrm{supp}(\pi)$, then by
Urysohn's Lemma there is a $\varphi\in\scc(M_2)$ so that
$\varphi=0$ on $\mathrm{supp}(\pi)$ and $\varphi(h_0)=1$. Since
$A\mapsto\hat{A}$ is an isomorphism, there is an $A\in\sca_2$ so
that $\hat{A}=\varphi$. Then by (\ref{plancherelmeasure})
$$\|A\delta_o\|_2^2=\lan
A,A\ran=\int_{\mathrm{supp}(\pi)}|\hat{A}(h)|^2d\pi(h)=0,$$ and
so $A=0$ by Lemma~\ref{inner}, contradicting $\hat{A}=\varphi\neq
0$.
\end{proof}

\section{Calculating the Plancherel Measure and the $\ell^2$-spectrum}\label{plan}

In this section we will calculate the Plancherel measure of
$\sca_2$. It turns out that there are two cases to consider. We
will then use these results to compute the $\ell^2$-spectrum of
$\sca$. The Plancherel measure will also be used in the proof of
Theorem~\ref{mainresult}, where we show that $h_u=h_{u}'$ for all
$u\in\Hom(P,\bc^{\times})$.

\begin{lem}\label{clas} $\tau_{\alpha}<1$ for some $\alpha\in
R$ if and only if $R=BC_n$ and $q_n<q_0$.
\end{lem}

\begin{proof} If $R$ is reduced we have $\tau_{\alpha}=q_{\alpha}$ for all $\alpha\in R$. Thus $\tau_{\alpha}<1$ for some
$\alpha\in R$ implies that $R=BC_n$ for some $n\geq1$. Thus $R$
may be described as follows (see \cite[VI, \ts4,
No.14]{bourbaki}). Let $E=\br^n$ with standard basis
$\{e_i\}_{i=1}^{n}$, and let $R$ consist of the vectors $\pm
e_i$, $\pm2e_i$ and $\pm e_j\pm e_k$ for $1\leq i\leq n$ and
$1\leq j<k\leq n$. Recall from \cite[Appendix]{p} that
in an affine building of type $BC_n$ we have
$q_1=\cdots=q_{n-1}$. Thus by the definition of the numbers
$\tau_{\alpha}$ we have $\tau_{\pm e_i}=q_nq_0^{-1}$,
$\tau_{\pm2e_i}=q_0$ and $\tau_{\pm e_j\pm e_k}=q_1$ for $1\leq
i\leq n$ and $1\leq j<k\leq n$. The result follows.
\end{proof}

Following \cite[Chapter~V]{macsph} we call the situation where
$\tau_{\alpha}\geq1$ for all $\alpha\in R$ the \textit{standard
case}, and we call the situation where $\tau_{\alpha}<1$ for some
$\alpha\in R$ the \textit{exceptional case}.

\subsection{The Standard Case} Let $\mathbb{U}=\{u\in\Hom(P,\bc^{\times}):|u^{\la}|=1\textrm{
for all $\la\in P$}\}$. Writing $u_i=u^{\la_i}$ for each
$i=1,\ldots,n$, we have $\mathbb{U}\cong\mathbb{T}^n$ where
$\bt=\{z\in\bc:|z|=1\}$.

In the next theorem we introduce (following \cite{macsph}) a
measure $\pi_0$ which we will shortly see is closely related to
the Plancherel measure $\pi$ (in the standard case). We will
write $\hat{A}(u)$ for $h_u(A)$ when $A\in\sca$ and $u\in\bu$. As
we shall see in Corollary~\ref{cor:cts}, each such $h_u$ is
continuous for the $\ell^2$--operator norm, and so
(\ref{eq:pi0st}) will also be valid for $A,B\in\sca_2$.

\begin{thm}\label{goodcase}(cf. \cite[Theorem~5.1.5]{macsph}) Let $du$ denote the normalised Haar measure on $\mathbb{U}$, and
let $\pi_0$ be the measure on $\mathbb{U}$ given by
$d\pi_0(u)=\frac{W_0(q^{-1})}{|W_0|}|c(u)|^{-2}du$. Then
\begin{align}\label{eq:pi0st}
\lan
A,B\ran=\int_{\mathbb{U}}\hat{A}(u)\overline{\hat{B}(u)}d\pi_0(u)\quad\textrm{for
all $A,B\in\sca$}.
\end{align}
\end{thm}

\begin{proof} We may assume that $A=A_{\mu}$ and $B=A_{\nu}$, where $\mu,\nu\in
P^+$. Then the integrand in (\ref{eq:pi0st}) is
$\widehat{A_{\mu}A_{\nu^*}}(u)$. Now
$A_{\mu}A_{\nu^*}=\sum_{\la\in P^+}a_{\mu,\nu^*;\la}A_{\la}$, and
since $a_{\mu,\nu^*;\la}=\delta_{\mu,\nu}/N_{\mu}$, it suffices
to show that
$\int_{\mathbb{U}}\widehat{A}_{\la}(u)d\pi_0(u)=\delta_{\la,0}$
for each $\la\in P^+$. Notice that if $u\in\mathbb{U}$, then
$\overline{u}=u^{-1}$, and so
$$|c(u)|^2=c(u)c(u^{-1})=\prod_{\alpha\in R}\frac{1-\taua^{-1}\tauah u^{-\cha}}{1-\tauah
u^{-\cha}}.$$ Thus $|c(wu)|^2=|c(u)|^2$ for all $w\in W_0$.
Furthermore, if $f(u)=\sum_{\la}a_{\la}u^{\la}$ is such that
$\sum_{\la}|a_{\la}|<\infty$, then $\int_{\mathbb{U}}f(u)du=a_0$.
It follows that
$\int_{\mathbb{U}}f(wu)du=\int_{\mathbb{U}}f(u)du$ for all $w\in
W_0$.

Using these facts we see that
\begin{align}
\label{standard}\int_{\mathbb{U}}\widehat{A}_{\la}(u)d\pi_0(u)&=q_{t_{\la}}^{-1/2}\int_{\mathbb{U}}\frac{u^{\la}}{c(u^{-1})}du.
\end{align}
Let $R_{\tau}^+=\{\alpha\in R^+\mid \tau_{\alpha}\neq 1\}$. Then
it is clear that we can write
\begin{align*}
\frac{1}{c(u^{-1})}=\prod_{\alpha\in R_{\tau}^+}\frac{1-\tauah
u^{\cha}}{1-\taua^{-1}\tauah u^{\cha}}=\sum_{\gamma\in
Q^+}a_{\gamma}u^{\gamma}
\end{align*}
where $a_{0}=1$ and the series is uniformly convergent. Since
$\{\la_i\}_{i=1}^n$ forms an \textit{acute} basis of $E$ \cite[VI,
\ts1, No.10]{bourbaki} we have $\lan\la,\la_i\ran\geq0$ for all
$\la\in P^+$ and for all $1\leq i\leq n$. Thus each $\la\in P^+$
is a linear combination of $\{\alpha_i\}_{i=1}^n$ with
nonnegative coefficients. It follows that if $\la\in P^+$,
$\gamma\in Q^+$ and $\la+\gamma=0$, then $\la=\gamma=0$. Hence by
(\ref{standard}) we have
$\int_{\mathbb{U}}\widehat{A}_{\la}(u)d\pi_0(u)=\delta_{\la,0}$,
completing the proof.
\end{proof}

Fix $o\in V_P$ and let $\ell^2_o(V_P)$ denote the space of all $
f\in\ell^2(V_P)$ which are constant on each set $V_{\la}(o)$. For
$A\in \sca_2$ define $\|A\|_o$ by
$$\|A\|_o=\mathrm{sup}\{\|Af\|_2\,:\,f\in\ell_o^2(V_P)\textrm{ and
}\|f\|_2\leq1\},$$ which defines a norm on $\sca_2$ (see
Lemma~\ref{inner}), and clearly $\|A\|_o\leq \|A\|$ for all
$A\in\sca_2$.

\begin{rem} In fact $\|A\|_o=\|A\|$ for all $A\in\sca_2$ (in both the standard and exceptional cases). To see this, recall that an injective homomorphism between two
$C^*$-algebras is an isometry \cite[Theorem~I.5.5]{davidson}. Let
$\Phi:\sca_2\to\mathscr{L}(\ell^2_o(V_P))$ be the linear map
given by $A\mapsto A|_{\ell^2_o(V_P)}$. Since
$\|A\|_o=\|A|_{\ell^2_o(V_P)}\|$, it suffices to show that $\Phi$
is an injection. This is clear from Lemma~\ref{inner}, for
$\Phi(A)=0$ implies that $A\delta_o=0$, and so $A=0$.
\end{rem}

\begin{cor}\label{cor:cts} Each $h_u$, $u\in \mathbb{U}$, is continuous for the
$\ell^2$-operator norm.
\end{cor}

\begin{proof} We show that $|h_u(A)|\leq\|A\|_o$ for all $A\in\sca$ and $u\in\mathbb{U}$.
Suppose that this condition fails for some $u_0\in\mathbb{U}$ and
$A\in\sca$. Then there exists $\delta>0$ so that
$|h_{u_0}(A)|>(1+\delta)\|A\|_o$. Since $h_u(A)$ is a Laurent
polynomial in $u_1,\ldots,u_n$ there exists a neighbourhood $\cn$
of $u_0$ in $\mathbb{U}$ such that $|h_u(A)|>(1+\delta)\|A\|_o$
for all $u\in \cn$. Let $\cn'$ denote the set of $u\in\mathbb{U}$
such that $|h_u(A)|>(1+\delta)\|A\|_o$, so $W_0\cn'=\cn'$. Let
$\mathbb{U}/W_0$ denote the set of $W_0$ orbits in $\mathbb{U}$.
It is compact Hausdorff with respect to the quotient topology, and
$$\scc(\mathbb{U}/W_0)\cong\{\varphi\in\scc(\mathbb{U})\mid\varphi(wu)=\varphi(u)\textrm{ for
all $w\in W_0$ and $u\in\mathbb{U}$}\}.$$ Now there exists
$\varphi\in\scc(\mathbb{U}/W_0)$ not identically~$0$ which is~$0$
outside~$\cn'$. By the Stone-Weierstrass Theorem, for any given
$\epsilon>0$ there exists $B\in\sca$ so that
$\|\hat{B}-\varphi\|_{\infty}<\epsilon$, and choosing $\epsilon$
suitably small we can ensure that
$$\int_{\cn'}|\hat{B}(u)|^2d\pi_0(u)\geq\frac{1}{1+\delta}\int_{\mathbb{U}}|\hat{B}(u)|^2d\pi_0(u)>0.$$
Thus by (\ref{eq:pi0st})
\begin{align*}
\|AB\delta_o\|_2^2&=\lan AB,AB\ran\\
&=\int_{\mathbb{U}}|\hat{A}(u)|^2|\hat{B}(u)|^2d\pi_0(u)\\
&\geq(1+\delta)^2\|A\|_o^2\int_{\cn'}|\hat{B}(u)|^2d\pi_0(u)\\
&\geq(1+\delta)\|A\|_o^2\|B\delta_o\|_2^2,
\end{align*}
and so $\|Af\|_2\geq\sqrt{1+\delta}\|A\|_o\|f\|_2$ for
$f=B\delta_o$, contrary to the definition of $\|A\|_o$.
\end{proof}

\begin{cor}\label{sp1} In the standard case, $M_2=\{\tilde{h}_u\mid u\in \bu\}$.
Moreover, the map $\varpi:u\mapsto\tilde{h}_u$ induces a
homeomorphism $\bu/W_0\to M_2$ (where $\bu$ is given the
Euclidean topology, $\bu/W_0$ is given the quotient topology, and
$M_2$ is given the w$^*$-topology), and the Plancherel measure
$\pi$ is the image of the measure $\pi_0$ of
Theorem~\ref{goodcase} under~$\varpi$.
\end{cor}
\begin{proof} The w$^*$-topology on~$M_2$, defined using the functionals
$h\mapsto h(A)$, $A\in\sca_2$, is compact, and so agrees with the
topology defined using only the functionals $h\mapsto h(A)$,
$A\in\sca$, since the latter is Hausdorff. Since each $h_u(A)$
($A\in\sca$ fixed) is a Laurent polynomial in $u_1,\ldots,u_n$,
the map $\varpi:u\mapsto\tilde{h}_u$, defined from~$\bu$ to~$M_2$
in light of Corollary~\ref{cor:cts}, is continuous. Thus
$\varpi(\bu)$ is closed in $M_2$, and $\varpi$ induces a
homeomorphism $\bu/W_0\to \varpi(\bu)$ since $\bu/W_0$ is
compact. The image~$\pi$ of $\pi_0$ under~$\varpi$ satisfies the
defining properties of the Plancherel measure. Since
$M_2=\mathrm{supp}(\pi)$ by Proposition~\ref{plancherelthm},
$$
\pi(M_2\backslash\varpi(\bu))=\pi_0(\varpi^{-1}(M_2\backslash\varpi(\bu)))=0,
$$
and so $M_2=\varpi(\bu)$. Thus $\varpi$ is surjective, and it is
injective by Proposition~\ref{inc}.
\end{proof}

\subsection{The Exceptional Case} Let $R=BC_n$ for some $n\geq1$ and suppose that $q_n<q_0$. Recall
the description of $R$ from the proof of Lemma~\ref{clas}. Let
$\alpha_i=e_i-e_{i+1}$ for $1\leq i\leq n-1$ and let
$\alpha_n=e_n$. The set $B=\{\alpha_i\}_{i=1}^{n}$ is a base of
$R$, and the set of positive roots with respect to $B$ is
$$R^+=\{e_i,2e_i,e_j-e_k,e_j+e_k\mid 1\leq i\leq n,\,1\leq j<k\leq
n\}.$$ Recall (from the proof of Lemma~\ref{clas}) that
$q_1=\cdots=q_{n-1}$ in this case. Let $a=\sqrt{q_nq_0}$ and
$b=\sqrt{q_n/q_0}$ (so $b<1$).

Let $u\in\Hom(P,\bc^{\times})$. Since $e_i\in P$ for each
$i=1,\ldots,n$, we may define numbers $t_i=t_i(u)$ by
$t_i=u^{e_i}$. We will now give a formula for $c(u)$ in this case
in terms of the numbers $\{t_{i}\}_{i=1}^n$ (see Remark~\ref{par}
for a related discussion).

If $\alpha=2e_i$, $1\leq i\leq n$, then $\alpha\in R_1\backslash
R_3$, and so $\tau_{\alpha}=q_0$. Now $\alpha/2=e_i\in
R_2\backslash R_3$, and so
$\tau_{\alpha/2}=q_{\alpha/2}q_0^{-1}$. Since
$|\alpha/2|=|\alpha_n|$ we have $q_{\alpha/2}=q_{\alpha_n}=q_n$,
and thus $\tau_{\alpha/2}=q_nq_0^{-1}$. Now if $\alpha=e_i$,
$1\leq i\leq n$, then $\alpha\in R_2\backslash R_3$, and so by
the above $\tau_{\alpha}=q_nq_0^{-1}$, and since $\alpha/2\notin
R$ we have $\tau_{\alpha/2}=1$. Since $(2e_i)^{\vee}=e_i$ and
$e_i^{\vee}=2e_i$, the factors in $c(u)$ (see (\ref{macsph}))
corresponding to the roots $\alpha=2e_i$ and $\alpha=e_i$ are
$$\frac{1-q_n^{-1/2}q_0^{-1/2}t_i^{-1}}{1-q_n^{-1/2}q_0^{1/2}t_i^{-1}}\cdot\frac{1-q_n^{-1}q_0t_i^{-2}}{1-t_i^{-2}}=\frac{(1-a^{-1}t_i^{-1})(1+b^{-1}t_i^{-1})}{1-t_i^{-2}}.$$
If $\alpha=e_j\pm e_k$, $1\leq j<k\leq n$, then $\alpha\in R_3$,
and so $\tau_{\alpha}=q_{\alpha}$. Since
$|\alpha|=|e_1-e_2|=|\alpha_1|$ we have
$q_{\alpha}=q_{\alpha_1}=q_1(=q_2=\cdots=q_{n-1})$, and so the
product of the two factors of $c(u)$ corresponding to the roots
$\alpha=e_j-e_k$ and $\alpha=e_j+e_k$ ($1\leq j<k\leq n$) is
$$\frac{(1-q_1^{-1}t_j^{-1}t_k)}{(1-t_j^{-1}t_k)}\cdot\frac{(1-q_1^{-1}t_j^{-1}t_k^{-1})}{(1-t_j^{-1}t_k^{-1})}.$$
Combining all these factors we see that $c(u)$ equals
$$\bigg\{\prod_{i=1}^{n}\frac{(1-a^{-1}t_i^{-1})(1+b^{-1}t_i^{-1})}{1-t_i^{-2}}\bigg\}\bigg\{\prod_{1\leq j<k\leq
n}\frac{(1-q_1^{-1}t_j^{-1}t_k)(1-q_1^{-1}t_j^{-1}t_k^{-1})}{(1-t_j^{-1}t_k)(1-t_j^{-1}t_k^{-1})}\bigg\}.$$

Notice that when $R=BC_n$, $Q=P$, and so we will be able to apply
the results of \cite{macsph} (see the paragraph after the proof
of Proposition~\ref{inc}). Thus when $q_1>1$ the Plancherel
measure here depends on how many of the numbers $q_1^kb$,
$k\in\bn$, are less than 1 (see \cite[page~70]{macsph}). Since we
have an underlying building we have the following simplification.
\begin{lem}\label{hig} If $q_1>1$, then $q_1b\geq 1$.
\end{lem}
\begin{proof} By a well known theorem of D. Higman (see \cite[page 30]{ronan} for
example), in a finite thick generalised $4$-gon with parameters
$(k,l)$, we have $k\leq l^2$ and $l\leq k^2$. Thus by
\cite[Theorem~3.5 and Proposition~3.2]{ronan} we have $q_1^2\geq
q_0$ (even if $q_0=1$), and so $q_1b\geq\sqrt{q_n}\geq1$.
\end{proof}

Let $\bt=\{z\in\bc:|z|=1\}$. Let $dt=dt_1\cdots dt_n$, where
$dt_i$ is normalised Haar measure on $\bt$. Define
$\phi_0(u)=c(u)c(u^{-1})$ and
$$\phi_1(u)=\lim_{t_1\to-b}\frac{\phi_0(u)}{1+b^{-1}t_1}\qquad\textrm{and}\qquad
dt'=d\delta_{-b}(t_1)dt_2\cdots dt_n.$$
Note that this limit exists since there is a factor $1+b^{-1}t_1$ in
$c(u^{-1})$. We use the isomorphism
$\mathbb{U}\to\bt^n$, $u\mapsto(t_1,\ldots,t_n)$ to identify
$\mathbb{U}$ with $\bt^n$. Define
$\mathbb{U}'=\{-b\}\times\bt^{n-1}$, and write
$U=\mathbb{U}\cup\mathbb{U}'$.

\begin{thm}\label{badcase} Let $\pi_0$ be the measure on $U=\mathbb{U}\cup\mathbb{U}'$ given by
$d\pi_0(u)=\frac{W_0(q^{-1})}{|W_0|}\frac{dt}{\phi_0(u)}$ on
$\mathbb{U}$ and
$d\pi_0(u)=\frac{W_0(q^{-1})}{|W_0'|}\frac{dt'}{\phi_1(u)}$ on
$\mathbb{U}'$, where $W_0'$ is the Coxeter group $C_{n-1}$ (with
$C_1=A_1$ and $C_0=\{1\}$). Then (in the exceptional case)
\begin{align}\label{eq:pi0ex}
\lan
A,B\ran=\int_{U}\hat{A}(u)\overline{\hat{B}(u)}d\pi_0(u)\quad\textrm{for
all $A,B\in\sca$}.
\end{align}
\end{thm}

\begin{proof} When $q_1>1$ this follows from the `group free' calculations made in \cite[Theorem~5.2.10]{macsph},
taking into account Lemma~\ref{hig}. If $q_1=1$ the formula for
$c(u)$ simplifies considerably, and a calculation similar to that
in \cite[Theorem~5.2.10]{macsph} proves the result in this case
too.
\end{proof}

As in the standard case we have the following corollary (see
Corollary~\ref{sp1}).

\begin{cor}\label{sp2} In the exceptional case, $M_2=\{\tilde{h}_u\mid u\in U\}$.
Moreover, the map $\varpi:u\mapsto\tilde{h}_u$ induces a
homeomorphism $(\bu/W_0)\cup(\bu'/W_0')\to M_2$ and the
Plancherel measure $\pi$ is the image of the measure $\pi_0$ of
Theorem~\ref{badcase} under~$\varpi$.
\end{cor}

\section{Equality of $h_u$ and $h_u'$}

In this section we show that $h_u=h_{u}'$ for all
$u\in\Hom(P,\bc^{\times})$, where $h_u'$ is as in
Theorem~\ref{main} (see also Corollary~\ref{2ndform}). To
conveniently state our results we will write $U=\mathbb{U}$ in
the standard case. Thus $U=\mathbb{U}$ in the standard case and
$U=\mathbb{U}\cup\mathbb{U}'$ in the exceptional case.

\begin{lem}\label{5.1} Let $\la,\mu,\nu\in P^+$. Then
\begin{enumerate}
\item[$\mathrm{(i)}$]$\int_{U}P_{\la}(u)\overline{P_{\mu}(u)}d\pi_0(u)=\delta_{\la,\mu}N_{\la}^{-1}$,
and
\item[$\mathrm{(ii)}$]$a_{\la,\mu;\nu}=N_{\nu}\int_{U}P_{\la}(u)P_{\mu}(u)\overline{P_{\nu}(u)}d\pi_0(u)$.
\end{enumerate}
\end{lem}

\begin{proof} (i) Since each $h_u$, $u\in U$, is continuous with respect to the $\ell^2$-operator norm
(see Corollary~\ref{sp1} and Corollary~\ref{sp2}) we have
$P_{\la}(u)=\widehat{A}_{\la}(u)$ for all $\la\in P^+$ and all
$u\in U$. Thus by (\ref{eq:pi0st}) in the standard case, and
(\ref{eq:pi0ex}) in the exceptional case,
$$\int_{U}P_{\la}(u)\overline{P_{\mu}(u)}d\pi_0(u)=\lan
A_{\la},A_{\mu}\ran=N_{\la}^{-1}\delta_{\la,\mu}.$$

(ii) Using the previous part we have
\begin{align*}
\int_{U}P_{\la}(u)P_{\mu}(u)\overline{P_{\nu}(u)}d\pi_0(u)&=\sum_{\eta\in
P^+}\left(a_{\la,\mu;\eta}\int_{U}P_{\eta}(u)\overline{P_{\nu}(u)}d\pi_0(u)\right)=N_{\nu}^{-1}a_{\la,\mu;\nu},
\end{align*}
completing the proof.
\end{proof}

\begin{thm}\label{mainresult} $h_u'=h^{\vphantom{'}}_u$ for all $u\in\Hom(P,\bc^{\times})$.
\end{thm}

\begin{proof} From (\ref{triang}) we have
\begin{align}\label{a2s}
h_{u}(A_{\la})=P_{\la}(u)=\sum_{\mu\in
\Pi_{\la}}a_{\la,\mu}u^{\mu}
\end{align}
for some numbers $a_{\la,\mu}$. On the other hand, by
Theorem~\ref{main} and Corollary~\ref{2ndform} we have
\begin{align}\label{eq:a2s}
h_u'(A_{\la})=\sum_{\mu\in\Pi_{\la}}r^{-\mu}a_{\la,\nu-\mu;\nu}u^{\mu}\qquad\textrm{for
any $\nu\ggg\la$}.
\end{align}
We will show that for all $\mu\in \Pi_{\la}$,
$a_{\la,\mu}=r^{-\mu}a_{\la,\nu-\mu;\nu}$ provided that
$\nu\ggg\la$. Comparing formulae (\ref{a2s}) and (\ref{eq:a2s})
this evidently proves that $h_u=h'_u$ for all
$u\in\Hom(P,\bc^{\times})$.

Let us first consider the standard case, so $U=\mathbb{U}$. Let
$\mu\in \Pi_{\la}$ and $\nu\ggg\la$. By Corollary~\ref{top} we
have $r^{-\mu}=\sqrt{N_{\nu-\mu}/N_{\nu}}$, and so by (\ref{ns})
we compute
$r^{-\mu}N_{\nu}=W_0(q^{-1})q_{t_{\nu}}^{1/2}q_{t_{\nu-\mu}}^{1/2}$.
Thus by Lemma~\ref{5.1}(ii),
\begin{align}\label{tttt}\begin{aligned}
r^{-\mu}a_{\la,\nu-\mu;\nu}&=r^{-\mu}N_{\nu}\int_{\mathbb{U}}P_{\la}(u)P_{\nu-\mu}(u)\overline{P_{\nu}(u)}d\pi_0(u)\\
&=W_0(q^{-1})q_{t_{\nu}}^{1/2}\int_{\mathbb{U}}P_{\la}(u)\overline{P_{\nu}(u)}\frac{u^{\nu-\mu}}{c(u^{-1})}du.
\end{aligned}\end{align} Since $|c(wu^{-1})|^2=|c(u^{-1})|^2$ for all
$u\in\mathbb{U}$ we see that $c(wu^{-1})/c(u^{-1})\in
L^1(\mathbb{U})$ for all $w\in W_0$, and so by (\ref{tttt}) we
have
\begin{align}\label{eq:refto}
r^{-\mu}a_{\la,\nu-\mu;\nu}&=\sum_{w\in
W_0}\int_{\mathbb{U}}u^{\nu-w\nu}\left(P_{\la}(u)u^{-\mu}\frac{c(wu^{-1})}{c(u^{-1})}\right)du.
\end{align}
We claim that the integral in (\ref{eq:refto}) is $0$ for all
$w\neq1$. To see this, notice that by Lemma~\ref{paper2lem4},
$r^{-\mu}a_{\la,\nu-\mu;\nu}$, $\mu\in\Pi_{\la}$, is independent
of $\nu\ggg\la$, and so we may choose $\nu=N(\la_1+\cdots+\la_n)$
for suitably large $N\in\bn$. Suppose $w\neq 1$. Since
$w(\la_1+\cdots+\la_n)\neq\la_1+\cdots+\la_n$ we see that
$$|\lan\la_1+\cdots+\la_n-w(\la_1+\cdots+\la_n),\alpha_{i_0}\ran|\geq1$$
for at least one $i_0\in I_0$, and so
$|\lan\nu-w\nu,\alpha_{i_0}\ran|\geq N$. Thus by the
Riemann-Lebesgue Lemma
$\lim_{N\to\infty}\int_{\mathbb{U}}u^{\nu-w\nu}f(u)du=0$ for all
$f\in L^1(\mathbb{U})$. Thus using (\ref{a2s}) we have
\begin{align*}
r^{-\mu}a_{\la,\nu-\mu;\nu}=\int_{\mathbb{U}}P_{\la}(u)u^{-\mu}du=a_{\la,\mu}.
\end{align*}

Let us now prove the result in the exceptional case, where
$U=\mathbb{U}\cup\mathbb{U}'$. Following the above we have
\begin{align}\label{eq:exi}
r^{-\mu}a_{\la,\nu-\mu;\nu}=a_{\la,\mu}+r^{-\mu}N_{\nu}\int_{\mathbb{U}'}P_{\la}(u)P_{\nu-\mu}(u)\overline{P_{\nu}(u)}d\pi_0(u)
\end{align}
whenever $\nu\ggg\la$. We will show that the integral in
(\ref{eq:exi}) is zero.

The group $W_0$ acts on $\{e_i\}_{i=1}^{n}$ as the group of all
signed permutations. Since $t_i=u^{e_i}$, for each $w\in W_0$ we
have
$w(t_1,\ldots,t_n)=(t_{\s_w(1)}^{\e_w(1)},\ldots,t_{\s_w(n)}^{\e_w(n)})$
where $\s_w$ is a permutation of $\{1,\ldots,n\}$ and
$\e_w:\{1,\ldots,n\}\to\{1,-1\}$. By directly examining the
formula for $c(u)$ we see that if $\e_w(\s_w(1))=-1$ then
$c(wu)|_{t_1=-b}=0$. Write $W_0^+=\{w\in W_0\mid
\e_w(\s_w(1))=1\}$. Note that for all $\la\in P^+$, $\la^*=\la$.
Thus $\overline{P_{\la}(u)}=P_{\la}(u)$ for all $u\in U$.
Following the calculation in the standard case, and using the
above observations, we see that
\begin{align*}
r^{-\mu}a_{\la,\nu-\mu;\nu}-a_{\la,\mu}&=r^{-\mu}N_{\nu}\int_{\mathbb{U}'}P_{\la}(u)P_{\nu-\mu}(u)P_{\nu}(u)d\pi_0(u)\\
&=\frac{1}{|W_0'|}\sum_{w,w'\in
W_0^+}\int_{\mathbb{U}'}u^{w\nu+w'\nu}\left(P_{\la}(u)u^{-w\mu}\frac{c(wu)c(w'u)}{\phi_1(u)}\right)dt'.
\end{align*}
Since $w_0=-1$ it is clear that if $w\in W_0^+$ then $w_0w\notin
W_0^+$. Take any pair $w,w'\in W_0^+$. As before, let
$\nu=N(\la_1+\cdots+\la_n)$ for sufficiently large $N$. Since
$w'\neq w_0w$ we have $w\nu+w'\nu\neq 0$. The same argument as in
the standard case now shows that $|\lan
w\nu+w'\nu,e_{i_0}\ran|\geq N$ for at least one $i_0\in I_0$.
Furthermore, since $w,w'\in W_0^+$ we have $\lan
w\nu+w'\nu,e_1\ran\geq0$. The result now follows by taking
$N\to\infty$, noting that $b<1$, and using the Riemann-Lebesgue
Lemma.
\end{proof}

Thus for all $\la\in P^+$ we have
\begin{align}\label{eq:twoforms}
h_{u}(A_{\la})=P_{\la}(u)=\int_{\O}(ur)^{h(x,y;\o)}d\nu_x(\o),\qquad
y\in V_{\la}(x).
\end{align}
As an application of (\ref{eq:twoforms}) we compute the norms
$\|A_{\la}\|$, $\la\in P^+$.

\begin{thm}\label{spectralradii} Let $\la\in P^+$. Then $\|A_{\la}\|=P_{\la}(1)$.
\end{thm}

\begin{proof} Since $A\mapsto\hat{A}$ is an isometry, by Corollaries~\ref{sp1}
and \ref{sp2} we have
\begin{align*}
\|A_{\la}\|=\|\hat{A}_{\la}\|_{\infty}=\sup\{|h(A_{\la})|:h\in
M_2\}=\sup\{|h_u(A_{\la})|:u\in U\},
\end{align*}
where, as usual, $U=\bu$ in the standard case and $U=\bu\cup\bu'$
in the exceptional case. In the standard case this implies that
$\|A_{\la}\|=P_{\la}(1)$, for by (\ref{eq:twoforms}) we have
$P_{\la}(1)>0$ and $|P_{\la}(u)|\leq P_{\la}(1)$ for all
$u\in\bu$ and $\la\in P^+$. In the exceptional case we only have
$\|A_{\la}\|\geq P_{\la}(1)$, and so it remains to show that
$\|A_{\la}\|\leq P_{\la}(1)$ in this case.

To see this, fix $o\in V_P$ and $\la\in P^+$. By
Theorem~\ref{main}(iii), Corollary~\ref{2ndform} and
(\ref{eq:twoforms}) we have
$(A_{\la}F_1^o)(x)=P_{\la}(1)F_{1}^o(x)$ for all $x\in V_P$.
Similarly, since $\la^*=\la$ here,
$$(A_{\la}^*F_1^o)(x)=(A_{\la^*}F_1^o)(x)=(A_{\la}F_1^o)(x)=P_{\la}(1)F_{1}^o(x)$$
for all $x\in V_P$. Since $F_1^o>0$ by (\ref{eq:twoforms}), the
Schur test (see \cite[p.103]{pedersen} for example) implies that
$\|A_{\la}\|\leq P_{\la}(1)$ (see also \cite[Lemma~4.1]{cm}).
\end{proof}

\begin{appendix}

\section{Calculation of $q_{t_{\la}}$}

The formula in Proposition~\ref{newlem} below is used in
Proposition~\ref{topt} to give an explicit formula for the
Radon-Nikodym derivative in Theorem~\ref{topology}. It is also
used after Proposition~\ref{inc} to compare our formula
(\ref{macsph}) with the formula in \cite[Theorem~4.1.2]{macsph}
(in the case $P=Q$).

\begin{prop}\label{newlem} Let $\la\in P^+$. Then
$$q_{t_{\la}}=\prod_{\alpha\in
R^+}\tau_{\alpha}^{\lan\la,\alpha\ran}.$$
\end{prop}

Before proving the above proposition, we need some preliminary
results.

\begin{lem}\label{Hs} Let $H$ be a wall of $\scx$. Suppose that $\pi_1$ is
a cotype $i$ (that is, cotype~$\{i\}$) panel of $H$ and that
$\pi_2$ is a cotype $j$ panel of $H$. Then $q_i=q_j$.
\end{lem}

\begin{proof} If $\s$ is any simplex of $\scx$ and $c$ any
chamber of $\scx$, then there is a unique chamber, denoted
$\mathrm{proj}_{\s}(c)$, nearest $c$ having $\s$ as a face
\cite[Corollary~3.9]{ronan}. We show that the map
$\varphi:\mathrm{st}(\pi_1)\to\mathrm{st}(\pi_2)$ given by
$\varphi(c)=\mathrm{proj}_{\pi_2}(c)$ is a bijection (here
$\mathrm{st}(\pi_i)$, $i=1,2$, denotes the set of chambers of
$\scx$ having $\pi_i$ as a face). Observe first that if
$c\in\mathrm{st}(\pi_1)$ then
$\mathrm{proj}_{\pi_1}(\mathrm{proj}_{\pi_2}(c))=c$. To see this,
let $\ca$ be any apartment containing $c$ and $H$ (see
\cite[Theorem~3.6]{ronan}), and let $H^+$ denote the half
apartment of $\ca$ containing $c$. Let $d$ be the unique chamber
in $\mathrm{st}(\pi_2)\cap H^+$. It follows from
\cite[Theorem~3.8]{ronan} that $\mathrm{proj}_{\pi_2}(c)=d$, and
so by symmetry $\mathrm{proj}_{\pi_1}(d)=c$. Similarly we have
$\mathrm{proj}_{\pi_2}(\mathrm{proj}_{\pi_1}(d))=d$ for all
$d\in\mathrm{st}(\pi_2)$. So the map $\varphi$ is bijective.
\end{proof}

Lemma~\ref{Hs} allows us to make the following (temporary)
definitions. Given a wall $H$ of $\scx$, write $q_H=q_i$, where
$i$ is the cotype of any panel of $H$. Now choose any apartment
$\ca$ of $\scx$, and let $\psi:\ca\to\S$ be a type rotating
isomorphism. For each $\alpha\in R$ and $k\in\bz$, write
$q_{\alpha;k}=q_H$, where $H=\psi^{-1}(H_{\alpha;k})$. We must
show that this definition is independent of the particular $\ca$
and $\psi$ chosen. To see this, let $\ca'$ be any (perhaps
different) apartment of $\scx$, and let $\psi':\ca'\to\S$ be a
type rotating isomorphism. Write $H'=\psi'^{-1}(H_{\alpha;k})$.
With $H$ as above, let $\pi$ be a panel of $H$, with cotype $i$,
say, and so $q_H=q_i$. The isomorphism
$\psi'^{-1}\circ\psi:\ca\to\ca'$ is type rotating and sends $H$
to~$H'$. Thus $(\psi'^{-1}\circ\psi)(\pi)$ is a panel of $H'$
with cotype $\s(i)$ for some $\s\in\Auttr(D)$, and so
$q_{H'}=q_{\s(i)}=q_i=q_H$, by \cite[Theorem~4.20]{p}.

\begin{lem}\label{lem:clas} Let $R$ be reduced. Then each wall of $\S$ contains an
element of $P$.
\end{lem}

\begin{proof} Each panel of each wall $H_{\alpha;k}$,
$\alpha\in R$, $k\in\bz$, contains $n-1$ vertices whose types are
pairwise distinct. Since $R$ is reduced, the good vertices are
simply the \textit{special vertices}, that is, the elements of
$P$ (see \cite[VI, \ts2, No.2, Proposition~3]{bourbaki}). Thus
when there are two or more good types the result follows.

This leaves the cases $E_8$, $F_4$ and $G_2$. Since
$H_{-\alpha;k}=H_{\alpha;-k}$ it suffices to prove the result
when $\alpha\in R^+$. Using the data in
\cite[Plates~VII-IX]{bourbaki} we see that for each
$\alpha=\sum_{i=1}^{n}a_i\alpha_i\in R^+$ there is an index $i_1$
such that $a_{i_1}=1$, or a pair of indices $(i_2,i_3)$ such that
$a_{i_2}=2$ and $a_{i_3}=3$. In the former case $k\la_{i_1}\in
H_{\alpha;k}$, and in the latter case $\frac{k}{2}\la_{i_2}\in
H_{\alpha;k}$ if $k$ is even, and
$\frac{k-3}{2}\la_{i_2}+\la_{i_3}\in H_{\alpha;k}$ if $k$ is odd.
\end{proof}

\begin{prop}\label{rel} If $R$ is reduced, then $q_{\alpha;k}=q_{{\alpha}}$ for all $\alpha\in R$ and
$k\in\bz$.
\end{prop}

\begin{proof} The proof consists of the following steps:
\begin{enumerate}
\item[(i)] $q_{w\alpha;0}=q_{\alpha;0}$ for all $\alpha\in R$ and
$w\in W_0$.
\item[(ii)] $q_{\alpha_i;0}=q_{\alpha_i}$ for each
$i=1,\ldots,n$.
\item[(iii)] $q_{\alpha;0}=q_{\alpha}$ for all $\alpha\in R$.
\item[(iv)] $q_{\alpha;k}=q_{\alpha;0}$ for all $\alpha\in R$ and
$k\in\bz$.
\end{enumerate}

(i) Let $\ca$ be an apartment of $\scx$, and let $\psi:\ca\to\S$
be a type rotating isomorphism. Write
$H=\psi^{-1}(H_{\alpha;0})$, so that $q_{\alpha;0}=q_H$. Let
$w\in W_0$. Now the isomorphism $\psi'=w\circ\psi:\ca\to\S$ is
type rotating, and
$\psi'^{-1}(H_{w\alpha;0})=\psi^{-1}(H_{\alpha;0})=H$. Thus
$q_{w\alpha;0}=q_H=q_{\alpha;0}$.

(ii) Let $C_0$ be the fundamental chamber of $\S$, and for each
$i=1,\ldots,n$ let $C_i=s_{i}C_0$. Let $\ca$ and $\psi:\ca\to\S$
be as in (i), and write $H=\psi^{-1}(H_{\alpha_i;0})$. Then
$\delta(\psi^{-1}(C_0),\psi^{-1}(C_i))=s_{\s(i)}$ for some
$\s\in\Auttr(D)$, and so $q_{\alpha_i;0}=q_H=q_{\s(i)}$, and so by
\cite[Theorem~4.20]{p} $q_{\alpha_i;0}=q_i=q_{\alpha_i}$.

(iii) Each $\alpha\in R$ is equal to $w\alpha_i$ for some $w\in
W_0$ and some $i$, and so (iii) follows from (i) and (ii).

(iv) Let $\alpha\in R$ and $k\in\bz$. By Lemma~\ref{lem:clas}
there exists $\la\in H_{\alpha;k}\cap P$, and so
$H_{\alpha;k}=t_{\la}(H_{\alpha;0})$. Let $\ca$ and $\psi$ be as
in (i), and write $H=\psi^{-1}(H_{\alpha;k})$, so that
$q_{\alpha;k}=q_{H}$. The map
$\psi'=t_{\la}^{-1}\circ\psi:\ca\to\S$ is a type rotating
isomorphism, and
$\psi'^{-1}(H_{\alpha;0})=\psi^{-1}(H_{\alpha;k})=H$. Thus
$q_{\alpha;k}=q_H=q_{\alpha;0}$.
\end{proof}

We need an analogue of Proposition~\ref{rel} when $R$ is of type
$BC_n$ for some $n\geq1$. Observe that if $\alpha\in
R_1\backslash R_3$, then $\alpha/2\in R_2\backslash R_3$, and
$H_{\alpha;2k}=H_{\alpha/2;k}$ for all $k\in\bz$. Thus we define
$\ch_1=\{H_{\alpha;k}\mid\alpha\in R_1\backslash R_3, \,k\textrm{
odd}\}$, $\ch_2=\{H_{\alpha;k}\mid \alpha\in R_2\backslash
R_3,\,k\in\bz\}$ and  $\ch_3=\{H_{\alpha;k}\mid\alpha\in R_3\}$.
Then $\ch=\ch_1\cup\ch_2\cup\ch_3$, where the union is disjoint.
We have
\begin{align}\label{eq:redthm}
q_{\alpha;k}=\begin{cases}q_0&\textrm{if $H_{\alpha;k}\in\ch_1$}\\
q_n&\textrm{if $H_{\alpha;k}\in\ch_2$}\\
q_{\alpha}&\textrm{if $H_{\alpha;k}\in\ch_3$}\end{cases}
\end{align}
We omit the details of this calculation.

\begin{rem} Proposition~\ref{rel} and formula (\ref{eq:redthm})
give the connection between our definitions of $R$ and
$\tau_{\alpha}$ and Macdonald's definitions of $\S_1$ and $q_a$
\cite[\ts3.1]{macsph}.
\end{rem}

Recall that for $\tilde{w}\in\tilde{W}$ we define
$\ch(\tilde{w})=\{H\in\ch\mid H\textrm{ separates $C_0$ and
$\tilde{w}C_0$}\}.$ Also, observe that each $H\in\ch$ is equal to
$H_{\alpha;k}$ for some $\alpha\in R_1^+$ and some $k\in\bz$, and
if $H_{\alpha;k}=H_{\alpha';k'}$ with $\alpha,\alpha'\in R_1^+$
and $k,k'\in\bz$, then $\alpha=\alpha'$ and $k=k'$.

\begin{proof}[Proof of Proposition~\ref{newlem}] Write
$t_{\la}=t_{\la}'g_l$, where $t_{\la}'\in W$ and $l=\tau(\la)$.
Then $\ch(t_{\la})=\ch(t_{\la}')$ and $q_{t_{\la}}=q_{t_{\la}'}$.
Suppose that $t_{\la}'=s_{i_1}\cdots s_{i_m}$ is a reduced
expression for $t_{\la}'$. Writing $H_{i}=H_{\alpha_i;0}$ if
$i=1,\ldots,n$ and $H_0=H_{\tilde{\alpha};1}$, we have
\begin{align}\label{eq:set}\begin{aligned}
\ch(t_{\la})&=\{H_{i_1},s_{i_1}H_{i_2},\ldots,s_{i_1}\cdots
s_{i_{m-1}}H_{i_m}\}\\
&=\{H_{\alpha;k_{\alpha}}\mid\alpha\in R_1^+\textrm{ and }0<
k_{\alpha}\leq\lan\la,\alpha\ran\}
\end{aligned}
\end{align}
where the hyperplanes in each set are pairwise distinct (see
\cite[Theorem~4.5]{h}). If $1\leq r\leq m$ and if $s_{i_1}\cdots
s_{i_{r-1}}H_{i_r}=H_{\alpha;k}$ then it is easy to see that
$q_{\alpha;k}=q_{i_r}$. Then using (\ref{eq:set}),
Proposition~\ref{rel}, and the fact that
$\lan\la,\alpha\ran\in2\bz$ for all $\alpha\in R_1\backslash
R_3$, we have
\begin{align*}
q_{t_{\la}}&=\prod_{r=1}^mq_{i_r}=\prod_{\alpha\in
R_1^+}\prod_{k_{\alpha}=1}^{\lan\la,\alpha\ran}q_{\alpha;k_{\alpha}}=\bigg[\prod_{\alpha\in
R_3^+}q_{\alpha}^{\lan\la,\alpha\ran}\bigg]\bigg[\prod_{\alpha\in
R_1^+\backslash R_3^+}(q_0q_n)^{\lan\la,\alpha\ran/2}\bigg],
\end{align*}
and the result follows by direct calculation.
\end{proof}

Proposition~\ref{equi} below is used to see that the formulae
(\ref{eq:macsphfn}) and \cite[(6.4)]{p} are equivalent. It also
gives the connection between our numbers $q_w$, $w\in W_0$, and
Macdonald's numbers $t_w$, $w\in W_0$, from~\cite[(3.8)]{loeth}
(which is needed in the proof of Lemma~\ref{llllls}).

Recall that $R_2=\{\alpha\in R\mid \frac{1}{2}\alpha\notin R\}$.
The \textit{inversion set} of $w\in W_0$ is
$$R_2(w)=\{\alpha\in R_2^+\mid\textrm{$H_{\alpha}$ is between
$C_0$ and $w^{-1}C_0$}\},$$ and we write $R(w)=\{\alpha\in
R^+\mid\textrm{$H_{\alpha}$ is between $C_0$ and $w^{-1}C_0$}\}.$

\begin{prop}\label{equi} For $w\in W_0$,
\begin{align}\label{forrr}
q_{w}=\prod_{\alpha\in R_2(w^{-1})}q_{\alpha}=\prod_{\alpha\in
R(w^{-1})}\tau_{\alpha}.
\end{align}
In particular,
$$q_{w_0}=\prod_{\alpha\in R^+}\tau_{\alpha}$$
(here $w_0$ is the longest element of $W_0$).
\end{prop}

\begin{proof} For $w\in W_0$ we have $R_2(w)=\{\alpha\in R_2^+\mid w\alpha\in
R^-\}$, and if $w=s_{i_1}\cdots s_{i_p}$ is a reduced expression
then
$$R_2(w)=\{\alpha_{i_p},s_{i_p}\alpha_{i_{p-1}},\ldots,s_{i_p}\cdots
s_{i_2}\alpha_{i_1}\}$$ (see \cite[VI, \ts1, No.6, Corollary~2 to
Proposition~17]{bourbaki}). The first equality in~(\ref{forrr})
follows.

The second equality in (\ref{forrr}) is clear if $R$ is reduced
(for $R_2=R$ and $\tau_{\alpha}=q_{\alpha}$). If $R$ is of type
$BC_n$, and if $\alpha\in R_2(w^{-1})\backslash R_3$, then
$q_{\alpha}=q_n=\tau_{\alpha}\tau_{2\alpha}$, verifying the
result in this case too.
\end{proof}

\section{Convex Hull}\label{B}

Given a subset $X\subset V_P$, define the \textit{convex hull of
$X$}, or $\mathrm{conv}(X)$, to be the set of good vertices that
lie in the intersection of all half-apartments that contain~$X$.

Let $H_{\alpha;k}$ be a hyperplane of $\S$. The (closed)
half-spaces of $\S$ associated to $H_{\alpha;k}$ are
$H_{\alpha;k}^+=\{z\in E\mid\lan z,\alpha\ran\geq k\}$ and
$H_{\alpha;k}^-=\{z\in E\mid\lan z,\alpha\ran\leq k\}$. This
gives an analogous definition of $\mathrm{conv}(X)$ for subsets
$X\subset P$, with the word half-apartment replaced by half-space.

Let $\leq$ denote the partial order on $P^+$ given by
$\mu\leq\la$ if and only if $\la-\mu\in P^+$. Note that this is
quite different to the partial order $\preceq$ on $P$ used
earlier.

\begin{lem}\label{first} Let $\la\in P^+$. Then
$\mathrm{conv}\{0,\la\}=\{\mu\in P^+\mid\mu\leq\la\}.$
\end{lem}

\begin{proof} Let $\mu\leq\la$ and write $\nu=\la-\mu\in
P^+$. Suppose that $0,\la\in H_{\alpha;k}^{\pm}$. Since
$H_{-\alpha;k}=H_{\alpha;-k}$ we may assume that $\alpha\in R^+$,
and since $0\in H_{\alpha;k}^{\pm}$ the only cases to consider are
$H_{\alpha;k}^-$ with $k\geq0$, and $H_{\alpha;k}^+$ with $k\leq
0$. In the case $k\geq0$ we have
$\lan\mu,\alpha\ran=\lan\la-\nu,\alpha\ran\leq\lan\la,\alpha\ran\leq
k$ and so $\mu\in H_{\alpha;k}^-$. In the case $k\leq 0$ we have
$\lan\mu,\alpha\ran\geq0\geq k$ and so $\mu\in H_{\alpha;k}^+$.
Thus $\{\mu\in
P^+\mid\mu\leq\la\}\subseteq\mathrm{conv}\{0,\la\}$.

Now suppose that $\mu\in\mathrm{conv}\{0,\la\}$. Since $0,\la\in
H_{\alpha_i;0}^+$ for each $i\in I_0$, we have $\mu\in
H_{\alpha_i;0}^+$ for each $i\in I_0$ too, and so $\mu\in P^+$.
Also, $0,\la\in H_{\alpha_i;\lan\la,\alpha_i\ran}^-$ for each
$i\in I_0$, and so $\lan\mu,\alpha_i\ran\leq\lan\la,\alpha_i\ran$
for each $i\in I_0$. That is, $\la-\mu\in P^+$, and so
$\mu\leq\la$.
\end{proof}

\begin{lem}\label{llllls} Let $\la,\mu\in P^+$. Then $|V_{\la}(x)\cap
V_{\mu^*}(y)|=1$ whenever $y\in V_{\la+\mu}(x)$.
\end{lem}

\begin{proof} We give a proof using calculations involving the
Macdonald spherical functions. We note that it is also possible
to give a purely `building theoretic' proof.

For $\la\in P^+$ let
$$P_{\la}'(x)=\frac{W_0(q^{-1})}{W_{0\la}(q^{-1})}q_{t_{\la}}^{1/2}P_{\la}(x)=q_{t_{\la}}^{-1/2}N_{\la}P_{\la}(x),$$
where we have used Proposition~\ref{betterf}, and define the
\textit{monomial symmetric function} $m_{\la}(x)\in\bc[P]^{W_0}$
by
$$m_{\la}(x)=\sum_{\mu\in W_0\la}x^{\mu}.$$
The set $\{m_{\la}(x)\}_{\la\in P^+}$ forms a basis for
$\bc[P]^{W_0}$, and by the calculations made in
\cite[\ts10]{loeth} (noting also Proposition~\ref{equi}) we have
\begin{align}\label{triang}
P'_{\la}(x)=\sum_{\mu\preceq\la}c_{\la,\mu}m_{\mu}(x),\,\textrm{
where } c_{\la,\la}=1.
\end{align}
This is the so called \textit{triangularity condition} of the
Macdonald spherical functions.

It is clear that
$m_{\la}(x)m_{\mu}(x)=\sum_{\nu\preceq\la+\mu}d_{\la,\mu;\nu}m_{\nu}(x)$
where $d_{\la,\mu;\la+\mu}=1$, and so it follows that
$$P_{\la}'(x)P_{\mu}'(x)=\sum_{\nu\preceq\la+\mu}e_{\la,\mu;\nu}P_{\nu}'(x)$$
for some numbers $e_{\la,\mu;\nu}$ with $e_{\la,\mu;\la+\mu}=1$.
It follows that for all $\la,\mu,\nu\in P^+$,
$$a_{\la,\mu;\nu}=q_{t_{\la}}^{1/2}q_{t_{\mu}}^{1/2}q_{t_{\nu}}^{-1/2}\frac{N_{\nu}}{N_{\la}N_{\mu}}e_{\la,\mu;\nu},$$
and since $q_{t_{\la+\mu}}=q_{t_{\la}}q_{t_{\mu}}$ we have
$a_{\la,\mu;\la+\mu}=N_{\la+\mu}N_{\la}^{-1}N_{\mu}^{-1}$. It
follows from (\ref{a}) that $|V_{\la}(x)\cap V_{\mu^*}(y)|=1$
whenever $y\in V_{\la+\mu}(x)$, proving the result.
\end{proof}

This immediately gives the following.

\begin{cor}\label{uni} Let $\la\in P^+$ and $\mu\leq\la$, and let $x,y\in V_P$ be any vertices with $y\in V_{\la}(x)$. There
exists a unique vertex, denoted $v_{\mu}(x,y)$, in the set
$V_{\mu}(x)\cap V_{\nu^*}(y)$, where $\nu=\la-\mu\in P^+$.
\end{cor}

\begin{thm}\label{ade} Let $\la\in P^+$, $x\in V_{P}$, and $y\in V_{\la}(x)$.
Then
$$\mathrm{conv}\{x,y\}=\{v_{\mu}(x,y)\mid \mu\leq\la\}.$$
\end{thm}

\begin{proof} Let $H$ be a half-apartment of $\scx$
containing $\{x,y\}$, and let $\ca$ be any apartment containing
$H$. It is easy to see (using Axiom~(B2) of \cite[p.76]{brown})
that there exists a type rotating isomorphism $\psi:\ca\to\S$
such that $\psi(x)=0$ and $\psi(y)=\la$. Let $\mu\leq\la$ and
write $\nu=\la-\mu\in P^+$. The vertex $v=\psi^{-1}(\mu)$ is in
both $V_{\mu}(x)$ and $V_{\nu^*}(y)$ (as $y\in V_{\nu}(v)$, for
$(t_{-\mu}\circ\psi)(v)=0$ and $(t_{-\mu}\circ\psi)(y)=\nu\in
P^+$), and so by Corollary~\ref{uni} $v_{\mu}(x,y)=v\in\ca$. Now
$\psi(H)$ is a half-space of $\S$ which contains $0$ and $\la$,
and so by Lemma~\ref{first} $\mu\in \psi(H)$. Thus
$v_{\mu}(x,y)=\psi^{-1}(\mu)\in H$, showing that
$$\{v_{\mu}(x,y)\mid \mu\leq\la\}\subseteq\mathrm{conv}\{x,y\}.$$

Suppose now that $v\in\mathrm{conv}\{x,y\}$. Thus there exists an
apartment $\ca$ containing $x,y$ and $v$. Let $\psi:\ca\to\S$ be a
type rotating isomorphism such that $\psi(x)=0$ and
$\psi(y)=\la$, and write $\mu=\psi(v)$. If
$\mu\notin\mathrm{conv}\{0,\la\}$ then there is a (closed)
half-space of $\S$ which contains $0$ and $\la$ but not $\mu$,
and it follows that there exists a half-apartment of $\ca$ which
contains $x$ and $y$ but not $v$, a contradiction. Thus
$\mu\in\mathrm{conv}\{0,\la\}$, and so by Lemma~\ref{first}
$\mu\leq\la$. It follows that $v\in V_{\mu}(x)\cap V_{\nu^*}(y)$,
where $\nu=\la-\mu\in P^+$, and so $v=v_{\mu}(x,y)$, completing
the proof.
\end{proof}

Note that the above shows that $\mathrm{conv}\{x,y\}$ is a finite
set for all $x,y\in V_P$. Indeed we have the following corollary,
which also shows that
$|\mathrm{conv}\{x,y\}|=|\mathrm{conv}\{u,v\}|$ whenever $y\in
V_{\la}(x)$ and $v\in V_{\la}(u)$, and that
$|\mathrm{conv}\{x,y\}|$ does not depend on the parameters of the
building.

\begin{cor} Let $y\in V_{\la}(x)$. Then
$|\{\mathrm{conv}\{x,y\}\}|=\prod_{i=1
}^n(\lan\la,\alpha_i\ran+1)$.
\end{cor}

\begin{proof} By Theorem~\ref{ade} we have
$|\{\mathrm{conv}\{x,y\}\}|=|\{\mu\in P^+\mid \mu\leq\la\}|$. Now
$\mu\leq\la$ if and only if $0\leq
\lan\mu,\alpha_i\ran\leq\lan\la,\alpha_i\ran$ for all
$i=1,\ldots,n$. Thus
$$|\{\mathrm{conv}\{x,y\}\}|=|\{(k_1,\ldots,k_n)\in\bn^n\mid 0\leq
k_i\leq\lan\la,\alpha_i\ran\textrm{ for all $i=1,\ldots,n$}\}|,$$
and the result follows.
\end{proof}

We conclude this appendix with a sketch of the following theorem,
which was used in the construction of the topology on $\O$.
Recall the definition of the maps $\varphi_{\mu,\la}$ made in the
paragraphs after the proof of Theorem~\ref{b}.

\begin{thm} Fix $x\in V_P$ and define
$\theta:\O\to\prod_{\la\in P^+}V_{\la}(x)$ by
$\o\mapsto(v_{\la}^x(\o))_{\la\in P^+}$. Then $\theta$ is a
bijection of $\O$ onto
$\varprojlim(V_{\la}(x),\varphi_{\mu,\la})$.
\end{thm}

\begin{proof} It is not too difficult to see that if $\cs$ and $\cs'$ are sectors with the same good vertices, then $\cs=\cs'$.
Thus it is clear that $\theta$ is injective. To show that $\t$ is
surjective, let $(v_{\nu})_{\nu\in
P^+}\in\varprojlim(V_{\la}(x),\varphi_{\mu,\la})$. For each
$m\geq1$ let $\mu_m=m(\la_1+\cdots+\la_n)$, and let $\cc(x;m)$
denote the set of chambers contained in the intersection of all
half-apartments containing $x$ and $v_{\mu_m}$. Since $\mu_m\in
P^{++}$ the sets $\cc(x;m)$ are nonempty for all $m\geq1$, and
$\cc(x;m)\subset \cc(x;k)$ whenever $m\leq k$. Furthermore, for
$m\geq1$ write $\cc_m$ for the set of chambers of $\S$ contained
in the intersection of all half-spaces containing $0$ and $\mu_m$.

For each $m\geq1$ there exists an apartment $\ca_m$ containing
$x$ and $v_{\mu_m}$, and a type rotating isomorphism
$\psi_m:\ca_m\to\S$ such that $\psi_m(x)=0$ and
$\psi_m(v_{\mu_m})=\mu_m$. Furthermore, if $\ca_m'$ and $\psi_m'$
also have these properties, then it is easy to see that
$\psi_m|_{\cc(x;m)}=\psi'_m|_{\cc(x;m)}$. Also,
$\psi_{m+1}|_{\cc(x;m)}=\psi_m|_{\cc(x;m)}$ for all $m\geq1$.

For each $m\geq1$ define $\xi_m:\cc_m\to\scx$ by
$\xi_m=\psi_m^{-1}|_{\cc_m}$. Since $\xi_{m+1}|_{\cc_m}=\xi_m$ we
have $\xi_{k}|_{\cc_m}=\xi_m$ for all $k\geq m$. We therefore
define $\xi:\cc(\cs_0)\to\scx$ by $\xi(C)=\xi_m(C)$ once
$C\in\cc_m$. By replacing the type map $\tau:V(\S)\to I$ on $\S$
by $\s_i\circ\tau$ where $i=\tau(x)$, we may take all of the above
isomorphisms to be type preserving, and so by
\cite[Theorem~3.6]{ronan} we see that $\xi$ extends to an isometry
$\tilde{\xi}:\cc(\S)\to\scx$. Then $\tilde{\xi}(\cc(\S))$ is an
apartment of $\scx$, and $\cs=\tilde{\xi}(\cc(\cs_0))$ is a
sector. Let $\o$ be the class of $\cs$. Then
$\t(\o)=(v_{\nu})_{\nu\in P^+}$.
\end{proof}

\end{appendix}

\bibliography{PhD.bib}
\bibliographystyle{plain}

\end{document}